\documentclass[12pt,a4paper]{amsart}
\usepackage{amsfonts}
\usepackage{amssymb}
\usepackage{bbm}
\usepackage{comment}
\usepackage{enumitem}
\usepackage{dsfont}
\usepackage{amsthm}
\usepackage{amsmath}
\usepackage{commath}
\usepackage{url}
\usepackage{epsfig}
\usepackage{subfigure}
\usepackage{adjustbox}
\usepackage{mathrsfs}  
\usepackage[ruled,vlined,linesnumbered]{algorithm2e}
\usepackage{graphicx}
\usepackage{a4wide}
\usepackage{todonotes}
\usepackage{mathtools}

\makeatletter
\@namedef{subjclassname@2010}{%
	\textup{2010} Mathematics Subject Classification}
\makeatother

\usepackage[T1]{fontenc}

\DeclareMathOperator{\N}{\mathbb{N}}

\newtheorem{thm}{Theorem}

\newtheorem{corollary}[thm]{Corollary}
\newtheorem{lemma}[thm]{Lemma}

\theoremstyle{definition}

\DeclareMathOperator{\desc}{desc}

\newcommand{\rn}{(r_n)_{n \in \mathbb{N}}}
\newcommand{\ls}{\leqslant}
\newcommand{\gs}{\geqslant}

\newcommand{\xn}{(x_n)_{n \in \mathbb{N}} }
\newcommand{\yn}{(y_n)_{n \in \mathbb{N}} }
\newcommand{\zn}{(z_n)_{n \in \mathbb{N}} }

\def \ba {\mathbf a}
\def \bb {\mathbf b}

\def \bN {\mathbb N}
\def \bP {\mathbb P}

\def \ba {\mathbf{a}}
\def \bb {\mathbf{b}}

\def \dd {\mathrm{d}}

\newcommand{\cG}{\mathcal{G}}

\newcommand{\cM}{\mathcal{M}}

\newcommand{\xin}{x_{i,N}}
\newcommand{\xiin}{x_{i+1,N}}

\newcommand{\yin}{y_{i,N}}
\newcommand{\yiin}{y_{i+1,N}}

\newcommand{\Ml}{\mathcal{M}_{\ell}^{(k)}}
\newcommand{\MlL}{\mathcal{M}_{\ell}^{(k)}(L)}
\newcommand{\MlR}{\mathcal{M}_{\ell}^{(k)}(R)}

\newcommand{\TlR}{\mathcal{T}_{\ell}^{(k)}(R)}
\newcommand{\TlL}{\mathcal{T}_{\ell}^{(k)}(L)}
\newcommand{\Tl}{\mathcal{T}_{\ell}^{(k)}}
\newcommand{\TjlR}{\mathcal{T}_{j(\ell)}^{(k)}(R)}

\newcommand{\Gl}{\mathcal{G}_{\ell}}

\newcommand{\nk}{(N_k)_{k\in \bN}}
\newcommand{\overbar}[1]{\mkern 2mu\widetilde{\mkern-3mu#1\mkern-1mu}\mkern 1mu}
\newcommand{\tMl}{\overbar{\mathcal{M}}_{\ell}^{(k)}}
\newcommand{\tMlL}{\tMl(L)}
\newcommand{\tMlR}{\tMl(R)}
\newcommand{\Mtl}{\mathcal{M}_{2\ell}^{(k)}}
\newcommand{\Mtjl}{\mathcal{M}_{j(\ell)}^{(k)}}
\newcommand{\tMtl}{\mathcal{N}_{2\ell}^{(k)}}
\newcommand{\tMtlo}{\mathcal{N}_{2\ell+1}^{(k)}}
\newcommand{\Mtlo}{\mathcal{M}_{2\ell+1}^{(k)}}

\newcommand{\prob}{\mathbb{P}}

\begin{document}
	\baselineskip=17pt
	\title{On Sequences with Exponentially Distributed Gaps }
	\author{Christoph Aistleitner}
	\address{TU Graz, Austria}
	\email{aistleitner@math.tugraz.at}
	\author{Manuel Hauke}
	\address{University of York, United Kingdom}
	\email{manuel.hauke@york.ac.uk}
	\author[A. Zafeiropoulos]{Agamemnon Zafeiropoulos}
	\address{Department of Mathematics, Technion, Haifa, Israel}
	\email{agamemn@campus.technion.ac.il}
	
	\date{}
	\thanks{ } 
	
	\begin{abstract} 
		It is well known that a sequence which has Poissonian correlations of all orders necessarily has exponentially distributed nearest-neighbor gaps. It is tempting to read this implication also in the other direction, by supposing that a sequence with exponential gap distribution must have Poissonian correlations, and as a consequence must be equidistributed. While the second of these assertions is known to be true (Poissonian correlations of any order indeed imply equidistribution), we show that the first assertion is generally false, by constructing a sequence which has exponential gap distribution but fails to be equidistributed (and as a consequence, also fails to have Poissonian correlations of any order). 
	\end{abstract}
	\maketitle
	
	\section{Introduction}
	Let $\xn \subseteq [0,1]$ be a sequence in the unit interval.  Given $N\gs 1$, we write  $x_{1,N} \ls x_{2,N} \ls \ldots \ls x_{N,N}$  for the points $x_1, \ldots, x_N$ written in increasing order. For convenience we also define $x_{N+1,N} = 1+x_{1,N}$, which can be interpreted as working on the unit torus instead of the unit interval. Then we call the numbers \[  \xiin - \xin \qquad (1\ls i \ls N)  \] the {\it gaps} of the sequence $\xn$ at stage $N$.  We say that the sequence $\xn$ has {\it exponential gap distribution} if 
	\[ \lim_{N\to \infty} \frac{1}{N}\#\Big\{ i\ls N : \xiin - \xin \ls \frac{s}{N} \Big\} = 1 - e^{-s} \quad \text{ for all } s>0. \]  Such behavior is called \emph{Poissonian}, in analogy with the gap structure of the Poisson process. A deterministic sequence in $[0,1]$ which has exponential gap  distribution can be seen as showing ``random'' behavior, since a sequence $(X_n)_{n\in\mathbb{N}}$ of independent random variables, all sampled from the uniform distribution on $[0,1]$, almost surely has exponentially distributed gaps. A standard argument to prove this assertion goes by showing that the correlations of all orders (see \cite{our_other_paper} for a definition) are almost surely Poissonian, and then expressing the gap distribution in terms of the correlation functions;  we refer to \cite[Appendix A]{kr} for a more detailed discussion. \par Due to its interpretation as a pseudorandomness property, the exponential gap distribution appears in numerous open problems. The famous {\it Berry--Tabor conjecture} \cite{berry} in mathematical physics asserts that the eigenvalues of (generic) completely integrable quantum dynamical systems have exponentially distributed gaps. In this direction, there exist several partial results, but as a whole, the conjecture is still widely open. Even in the metric sense the theory of asymptotic gap distributions is far from being understood, even though there has been some progress. For example, for the sequence $(\{n^2\alpha\})_{n\in\mathbb{N}}$ Rudnick and Sarnak \cite{rud_sar} proved Poissonian pair correlation for almost all $\alpha$, and alternative proofs are known \cite{heath,mark_stro}. Whether this sequence has Poissonian triple correlation for almost all $\alpha$ is still unknown, and seems to be a difficult problem (cf.\ \cite{tech_walk}). Rudnick, Sarnak and Zaharescu \cite{rudnick2} conjectured that for almost all $\alpha$ the sequence has Poissonian correlations of all orders, as well as exponentially distributed gaps (see also \cite[Conjecture 1.1]{marklof2}), but this is widely out of reach with current methods. \par
	There exist a few examples of sequences which are known to have exponentially distributed gaps. Rudnick, Sarnak and Zaharescu \cite{rudnick2} gave a Diophantine criterion for $\alpha$ which guarantees that the sequence $(\{n^2 \alpha\})_{n\in\mathbb{N}}$ has exponential gaps (without providing an explicit example of a number $\alpha$ satisfying the condition), and Lutsko and Technau proved that the sequence $(\{\alpha (\log n)^A\})_{n\in\bN}$ for $\alpha >0$ and $A>1$ has exponential gaps \cite{lutsko}. Among metric results, it is known that the sequence $(\{q_n \alpha\})_{n\in \mathbb{N}}$  for a lacunary $(q_n)_{n\in\mathbb{N}}$ has exponential gap distribution for almost all $\alpha$ \cite{cy,rudnick3}, and that $(\{\alpha^n\})_{n \in \mathbb{N}}$ has exponential gaps for almost all $\alpha>1$ \cite{abty}. \par
	The standard method of proof in all of these examples is the one described earlier for the random case, by showing that the correlations of all orders are Poissonian.  As a consequence, all sequences that are known to have exponential gaps are also uniformly distributed sequences, because Poissonian correlations of any fixed order $k\gs 2$ are known to imply equidistribution (also called uniform distribution mod $1$); see \cite{alp,sigrid,marklof,steinerberger}.  Thus it is natural to think of exponential gap distribution as being a stronger property than uniform distribution, and sometimes it seems to have been considered as a fact that this actually is the case (see e.g. \cite{lutsko}). It is easy to construct a sequence which is uniformly distributed but does not have exponential gap distribution. However, the question whether there can be a sequence which has exponential gap distribution but is not uniformly distributed modulo 1 was unclear until now. This is the main topic of the present paper. Among earlier papers devoted (at least in parts) to the connection between gap structure and correlation functions, we mention \cite{abr, allp, ac, larch, weiss}.\\
	
	\noindent Before we proceed to more details, we examine what happens when this question is posed for triangular arrays instead of sequences. A {\it triangular array} in $[0,1]$ is a sequence $(X_N)_{N\gs 1}$ where each $X_N := (x_{1}^N,\ldots,x_N^{N})$ is an $N$-tuple of elements in $[0,1]$. Clearly, any sequence $\xn$ can be interpreted as a special case of a triangular array by defining $X_N = (x_1,\ldots,x_N),$ but not all triangular arrays correspond to a sequence in such a way. \par The definitions of all well-known statistical properties of sequences can be extended to triangular arrays in a natural way. For example, we say that the triangular array $(X_N)_{N \gs 1}$ is uniformly distributed if for all $0 \ls a \ls b \ls 1$, we have  \[\lim_{N \to \infty} \frac{1}{N}\#\{n \ls N: x_{n}^N \in [a,b]\} = b-a, \]
	and one can define Poissonian pair correlations for triangular arrays in a similar fashion.  \par Marklof proved  that any triangular array $(X_N)_{N \gs 1}$ that has Poissonian pair correlations is uniformly distributed \cite{marklof}. By adapting the arguments of \cite{our_other_paper} it can easily be shown that  triangular arrays with Poissonian correlations of any order are also necessarily uniformly distributed. \par However, when it comes to the gap distribution of triangular arrays, it can be shown rather easily that exponential gaps do not necessarily imply uniform distribution. To see this, let $\xn$ be a sequence with  exponential gaps and write $\gamma_{1,N} \ls \ldots \ls \gamma_{N,N}$ for the gaps of this sequence at stage $N$, sorted in increasing order. Now define  a triangular array $(X_N)_{N \gs 1}$ by setting  $x_n^N = \sum_{j \ls n} \gamma_{j,N}, 1 \ls n \ls N$ for each $N\gs 1.$ Then $(X_N)_{N \gs 1}$ has exponential gaps because its gaps at each stage coincide with those of  $\xn.$  However, since the gaps appear in increasing order, one can show that for any $s>0,$ \[ \lim_{N\to\infty}\frac{1}{N}\#\left\{ n\ls N : x_n^N \ls F(s)\right\} = 1-e^{-s}, \quad \text{ where } F(s) = 1 - (1+s)e^{-s},  \]  whence $(X_N)_{N\gs 1}$ is not uniformly distributed mod $1$. This is because  moving the smaller gaps towards the ``left'' side of $[0,1]$ brought too many points close to $0,$ thereby destroying the property of uniform distribution.\par 	Clearly, this argument cannot be applied for sequences instead of triangular arrays, as the triangular array defined above does not correspond to a sequence.  Generally speaking, constructing an infinite sequence in $[0,1]$ with a desired asymptotic gap distribution is a very complicated problem (much more so than the corresponding problem for triangular arrays), since a new gap can only be created by splitting an existing gap into two smaller ones, thus actually creating two new gaps and removing an existing one. Furthermore, the scaling factor $1/N$ changes permanently as the total number of points is increased, so that existing gaps that remain untouched appear to grow in size when the total number of points is increased. Thus the question arises whether for sequences a certain asymptotic gap distribution determines the distribution of the sequence itself.
	\par The following theorem shows that indeed,  under some rather mild assumptions, the gap structure determines the distribution. In the formulation of the theorem, we say that a function $g:[0,1]\to [0,1]$ is  the {\it asymptotic density function} of a sequence $\xn\subseteq [0,1]$ if \[ \lim_{N\to\infty}\frac{1}{N}\#\{ n\ls N : x_n\ls x \} = \int_0^x g(t)\,\mathrm{d}t, \qquad (0\ls x\ls 1).  \]
	If $\xn$ and $\yn$ are two sequences with gaps $g_{i,N}=\xiin-\xin$ and $h_{i,N}=\yiin-\yin,$ respectively, we say that they have the same gaps if $\{g_{i,N}\}_{i\ls N} = \{h_{i,N}\}_{i \ls N}$ for all $N\gs 1.$ (Throughout the paper, the sets of gaps are interpreted as multisets, so equality means that the gap lengths appear in the same multiplicity.)  Furthermore, we say that the gaps $\cG = (\{g_{i,N}: i\ls N\})_{N\gs 1}$ of a sequence $\xn$ are {\it distinct} if for any $N\gs 1$ we have $g_{i,N} \neq g_{j,N}$ for $i\neq j.$ 
	
	\begin{thm}\label{gaps_sets_det_distr}
		Let $\xn \subseteq [0,1]$ be a dense sequence with distinct gaps $\cG = (\{g_{i,N}: i\ls N\})_{N\gs 1}.$ Further, assume that $\xn$ has an asymptotic density function $g$ such that $\int_{0}^1 g^k(x) \,\mathrm{d}x < \infty$ for all $k \in \mathbb{N}$. Let $\yn$ be a sequence with the same gaps as $\xn$. Then $\yn$ has an asymptotic density function $f$ which satisfies
		\begin{equation}\label{same_moments} \int_0^1 f(x)^k \,\mathrm{d} x = \int_0^1 g(x)^k \,\mathrm{d} x \qquad (k\gs 1).  \end{equation}
	\end{thm}
	Applying Theorem \ref{gaps_sets_det_distr} to a sequence $\xn$ that is uniformly distributed and has distinct gaps, we see that any sequence $\yn$ with the same gaps is automatically uniformly distributed as well. It might be tempting to think that the assumption of ``same gap structure and distinct gaps'' in Theorem \ref{gaps_sets_det_distr} can be replaced by ``same gap distribution'' alone, so that having the same gap distribution already implies \eqref{same_moments}. Since a random sequence has exponential gaps and is equidistributed almost surely, in view of \eqref{same_moments} one might be led to expect that any other sequence having exponential gaps must be equidistributed as well.
	
	However, somewhat surprisingly, it turns out that the situation is not quite so easy. In the main result of this paper, in Theorem \ref{thm1} below, we construct a sequence which has exponentially distributed gaps but fails to be equidistributed. The sequence is constructed in such a way that there are (many) gaps of exactly the same size, and this plays a crucial role in the argument. Since a new gap can only arise from splitting an existing gap into two parts, the assumption of having distinct gaps determines the structure of the sequence in a strong way. On the contrary, when gaps of exactly the same size are allowed, and if the sequence is constructed in such a way that it has many gaps of the same size, then for the creation of a new gap there is some freedom of choice about the location of the existing gap that is split into two. This ``choice of location'' (resulting from the existence of many gaps of the same size) is the key feature in our construction which allows us to avoid equidistribution of the resulting sequence.  
	
	\begin{thm}\label{thm1} There exists a sequence $\xn \subseteq [0,1]$  with exponentially distributed gaps that is not   uniformly distributed. \end{thm}
	\noindent  
	In fact, in order to prove Theorem \ref{thm1}, we construct two sequences $\xn,\zn$ which do not only have the same asymptotic gap distribution, but precisely the same gap structure. However, while $\zn$ is uniformly distributed, $\xn$ is not -- we note in passing that this implies that the assumption of having distinct gaps in the statement of Theorem \ref{gaps_sets_det_distr} is indeed necessary. Theorem \ref{thm1} refutes the perception that having exponential gap distribution is a ``stronger'' property than uniform distribution. In fact, as the theorem shows, the two properties are logically independent and all possible combinations of their logical values can be realized.\\
	
	As mentioned earlier, any sequence with Poissonian correlations of some given order $k\gs 2$ is uniformly distributed. We thus obtain the following statement as an immediate corollary.
	
	\begin{corollary}
		There exists a sequence $\xn \subseteq [0,1]$  with exponentially distributed gaps that does not have Poissonian correlations of any order $k\gs 2.$
	\end{corollary}
	
	\noindent {\bf Notation: } Throughout the text, we use the standard O- and o- notations. Given a set $A$, we shall use both $\# A$ and $|A|$ to denote its cardinality. Also given two random variables $X, Y$ on a common probability space, we shall write $X\stackrel{d}{=} Y$ to state that the distributions of $X$ and $Y$ coincide. \\
	
	\noindent {\bf Acknowledgments:} The first author was supported by the Austrian Science Fund (FWF), projects I-4945, I-5554, P-34763 and P-35322. The second author was supported by FWF project P-35322 and by the EPSRC grant EP/X030784/1. The third author was supported by European Research Council (ERC) under the European Union’s Horizon 2020 Research and Innovation Program, Grant agreement no. 754475. We would like to thank Lorenz Fr\"uhwirth for many valuable comments.
	
	\section{Proof of Theorem \ref{thm1}} \label{section_2}

	Before starting with the proof of Theorem \ref{thm1}, we give a brief heuristic description of the construction of the sequence $\xn$. We emphasize again that the key difficulty is to construct an \emph{infinite sequence} which has a prescribed gap statistic, and not just a finite point set (or a sequence of finite point sets, i.e.\ a triangular array). We start with a sequence $\yn$ of random numbers in $[0,1]$, which are picked independently with respect to the uniform distribution.  We will also need a (fast growing) sequence $(N_k)_{k\gs 1}$ of positive integers that create different ``stages'' of the construction of $\xn$.\par The  $k$--th stage involves the definition of the terms $x_n$ for $2N_{k-1} < n \ls 2N_k$. For $2N_{k-1} < n \ls N_k,$ the value of $x_n$ is that of $y_n$, rounded down at some suitable scale which depends on $k.$ The subsequent points  $x_n,\,\, N_k < n \ls 2N_k$, are defined by shifting the values of $y_n,\,\, N_k < n \ls 2N_k$, in a way such that the gap structure of $(y_n)_{n \in \N}$ is fully preserved, but such that there are (significantly) more points in the left half of the unit interval than in the right half. \par Since at each stage a sufficiently large proportion of points is shifted to the left half of the unit interval, the resulting  sequence will not be uniformly distributed modulo $1$. At the same time, rounding down and shifting the elements of $\yn$ is done in such a way so that  $\xn$ has asymptotically the same gap distribution as $\yn$, which is exponential due to the randomness of $\yn$. The rounding-down procedure at each stage is a key ingredient, which creates a large proportion of gaps of identical size (which allows us to shift some of the points without destroying the overall gap statistic). Note that as a consequence of this rounding-down process,  Theorem \ref{gaps_sets_det_distr} does not apply to the sequence $\xn$ constructed in this way.\\
	
	Now we give the details of the proof of Theorem \ref{thm1}. We begin with a sequence  $\yn$ of numbers that are picked independently and uniformly at random from the interval $[0,1]$. Formally speaking, we consider a sequence $\yn$  of independent random variables on a probability space $(\Omega, \Sigma, \bP ),$  all following the  uniform distribution  $U[0,1]$, and we work with ``random''
	sequences  $(y_n(\omega))_{n\in\mathbb{N}} \subseteq [0,1]$  which depend on $\omega\in \Omega$. Throughout the paper we will suppress this dependence of $\yn$ and all derived quantities on $\omega \in \Omega$, for the sake of readability. \par 
	Let $\nk$ be a strictly increasing (deterministic) sequence of positive integers such that  
	\begin{equation} \label{ratio}
		\begin{gathered}
			N_1 \text{ is even,} \quad  N_k \text{ divides } N_{k+1} \text{ for all } k\gs 1, \\ \text{ and } \quad  N_{k} = o( N_{k+1} ) \text{ as } k\to \infty. 
		\end{gathered}
	\end{equation}
	The  sequence $\xn$ will be defined recursively, with the $k$--th stage of the construction involving the definition of the terms $x_n$ for $2N_{k-1} < n \ls 2N_k.$  The numbers $\xn$ will depend on $\yn$, hence they are also random variables on $(\Omega, \Sigma, \bP )$.  We shall first present the definition of $\xn$, assuming that $\nk$ is given. We then prove that the sequence $\nk$ can be chosen in such a way that \eqref{ratio} is satisfied, and such that with positive probability, the terms of $\xn$ satisfy the conclusion of Theorem \ref{thm1}.\par   Before we proceed to the proof, we need to introduce one last definition. In what follows, for a given sequence $\xn$ we will need to associate the order points $x_{i,N}, ~1 \ls i \ls N,$ with the original (unordered) points $x_n, ~1 \ls n \ls N$. Given an integer $N\gs 1,$ we say that the point $x_{i,N}$ {\it comes from} $x_n \, (i, n\ls N)$ if $x_{i,N} = x_n.$ When several points coincide, i.e.\ when $x_{n_1} = x_{n_2} = \ldots = x_{n_p} = x_{i+1,N} = x_{i+2,N} = \ldots = x_{i+p,N}$ for some integers $1 \ls n_1 < \ldots < n_{n_p} \ls N$ and $p\ls N, i\ls N-p,$ we say that $x_{i+k,N}$ comes from the point $x_{n_k}, k=1,\ldots, p.$  \\

	\noindent {\bf The First  Step.} 
	For $n=1,2,\ldots, N_1$ set 
	\[   \tilde{y}_n = \frac{1}{10N_1}\lfloor 10 N_1 y_n \rfloor \qquad \text{ and } \qquad x_n = \tilde{y}_n. \]
	In other words, $x_n$ is the approximation to $y_n$ up to $10N_1$ decimal digits. By definition, the gaps between the points $x_n$ can only attain values \[  \frac{\ell}{10N_1}, \qquad 
	\ell \in \{0,\ldots,10N_1\}.\]
	For any $1\ls \ell \ls 10N_1$  consider the sets 
	\begin{align*}
		\cM_{\ell}^{(1)} &= \left\{ i\ls N_1 : x_{i+1,N_1} - x_{i,N_1} =\ell/(10N_1)   \right\}, \\
		\cM_{\ell}^{(1)}(L) &= \left\{ i \in \cM_{\ell}^{(1)} : i < N_k/2  \right\}, \\
		\cM_{\ell}^{(1)}(R) &= \left\{ i \in \cM_{\ell}^{(1)} : i \gs N_k/2   \right\}.
	\end{align*}
	We are now going to define the next points $x_{N_1+1}, \ldots, x_{2N_1}$ by suitably shifting some of the points $y_{n},\,\, N_1< n \ls 2N_1$.\par For each $1\ls \ell \ls 10N_1$, write $ \{ i_1, \ldots, i_{\mu_1(\ell)} \} $ and $\{ j_1, \ldots, j_{\nu_1(\ell)} \}$ for the sets 
	\begin{eqnarray*}
		\{  i\in \cM_{\ell}^{(1)}(L) & : & \#\left\{  N_1 < n \ls 2N_1 : y_n \in [x_{i,N_1}, x_{i+1,N_1})  \} =0  \right\} \quad \text{ and }\\
		\{  i\in \cM_{\ell}^{(1)}(R) &:&   \#\left\{  N_1 < n \ls 2N_1 : y_n \in [x_{i,N_1}, x_{i+1,N_1})  \right\} >0  \},
	\end{eqnarray*}
	respectively. Let $f_\ell : [0,1) \rightarrow [0,1)$ be the piecewise linear transformation defined as follows: if $\rho_1(\ell) = \min\{ \mu_1(\ell), \nu_1(\ell)\},$ then $f_{\ell}$ is the transformation that maps 
	\[  [x_{i_s, N_1}, x_{i_s+1,N_1} ) \quad \text{to} \quad  [x_{j_s, N_1}, x_{j_s+1,N_1} ) \qquad (s=1,\ldots, \rho_\ell)  \]
	and vice versa, and leaves the remaining points fixed. Note that by construction the two intervals in the formula above have the same length.
	
	We shall apply all such transformations $f_\ell,\,\, \ell =1,\ldots, 10 N_1$ successively: let \[ f^{(1)} = f_1 \circ \ldots \circ f_{10N_1}  \] and define  
	\[ \tilde{y}_n = \frac{1}{10N_1}\lfloor 10N_1 y_n \rfloor  \quad \text{ and } \quad  x_n = f^{(1)}(\tilde{y}_n) \qquad (N_1 < n \ls 2N_1).   \]
	Note how, heuristically, $f^{(1)}$ considers the intervals that arise from the points with index at most $N_1$, and then exchanges an interval in the \emph{left} half of $[0,1]$ which contains \emph{no} point with index in $(N_1, 2N_1]$ against an equal-length interval in the \emph{right} half of $[0,1]$ that \emph{does} contain one or more points with index in $(N_1,2N_1]$; thus, $f^{(1)}$ effectively shifts points from the right half of $[0,1]$ into the left half without changing the overall gap statistic. Note also how this procedure requires the existence of gaps of the same size, which can be exchanged against each other; in order to have a significant overall effect, leading to a non-uniform distribution of the resulting sequence $(x_n)_{n \in \mathbb{N}}$, we need to ensure that there are ``many'' gaps of equal size to which the procedure can be applied. \\
	
	\noindent {\bf The inductive step.} Assume that for some $k\gs 1$ the points $x_1, \ldots, x_{2N_{k-1}}$ have been defined. Given $N_k$ such that \eqref{ratio} holds, we set
	\[  \tilde{y}_n = \frac{1}{10^k N_k}\lfloor 10^k N_k y_n \rfloor \qquad \text{ and } \qquad x_n =\tilde{y}_n \qquad  (2N_{k-1} < n \ls N_{k}).   \]
	Similarly to the first step, we define the sets 
	\begin{align*}
		\cM_{\ell}^{(k)} \,\, &= \left\{  i\ls N_{k} : x_{i+1,N_k} - x_{i,N_k} =\ell/(10^kN_k)   \right\}, \\
		\cM_{\ell}^{(k)}(L) &= \left\{ i \in \cM_{\ell}^{(k)} : i < N_k/2   \right\}, \\
		\cM_{\ell}^{(k)}(R) &= \left\{ i \in \cM_{\ell}^{(k)} : i \gs N_k/2  \right\},
	\end{align*}
	where $ 0 \ls \ell \ls 10^k N_k$. Write $ \{ i_1, \ldots, i_\mu \} $ and $\{ j_1, \ldots, j_\nu \}$ for the sets 
	\begin{eqnarray*}
		\{  i\in \cM_{\ell}^{(k)}(L) & : & \#\{  N_k < n \ls 2N_k : y_n \in [x_{i,N_k}, x_{i+1,N_k})  \} =0  \} \quad \text{ and }\\
		\{  i\in \cM_{\ell}^{(k)}(R) &:&   \#\{  N_k < n \ls 2N_k : y_n \in [x_{i,N_k}, x_{i+1,N_k})  \} >0  \},
	\end{eqnarray*}
	respectively, where $\mu=\mu_k(\ell)$ and $\nu=\nu_k(\ell).$ Set 
	$\rho = \rho_k(\ell) = \min\{ \mu_k(\ell), \nu_k(\ell) \}$. We will define the transformation $f^{(k)}$ similar to the definition of $f^{(1)}$ in the first step of the construction. For each $\ell =1, \ldots, 10^k N_k$ we write $f_\ell^{(k)}: [0,1) \to [0,1)$ for the function that maps 
	\[  [x_{i_s,N_k}, x_{i_s+1,N_k}) \quad \text{ to } \quad [x_{j_s,N_k}, x_{j_s+1,N_k}) \qquad  (s=1,\ldots, \rho_k(\ell)) \]
	and vice versa and keeps the remaining points fixed. The map $f^{(k)}$ is defined as 
	\[ f^{(k)} = f_1^{(k)} \circ \ldots \circ f_{10^k N_k}^{(k)}.   \]
	We now set 
	\[ \tilde{y}_n = \frac{1}{10^k N_k}\lfloor 10^k N_k y_n \rfloor    \quad \text{ and }\quad x_n = f^{(k)}(\tilde{y}_n)  \qquad (N_{k} < n \ls 2N_k).     \]
	This completes the definition of the sequence $\xn.$

	\subsection{Exponential Gaps} In this part we prove that the sequence $\xn$ defined above has exponential gaps almost surely. 
	
	We first prove a closely related statement in a more abstract setup. Let $(N_k)_{k \in \N}$ be a sequence that satisfies \eqref{ratio}.  We consider two sequences $(Q_n)_{n\in \bN}$ and $\rn$ of random variables, such that all of them are mutually independent.  The distribution of the random variables is as follows:  for $1\ls n\ls N_1$, 
	\[ \begin{split} \bP(a<Q_n<b) &= (b-a)10N_1\qquad \left(0\ls a < b \ls \frac{1}{10N_1} \right) \quad \text{ and } \\  \bP\Big(r_n = \frac{i}{10N_1}\Big) &= \frac{1}{10N_1} \qquad  \left(0 \ls i < 10N_1 \right).\end{split} \]
	For $n>N_1,$ if $k\gs 2$ is the unique integer such that $2N_{k-1} < n \ls 2N_k,$ then 
	\[ \begin{split} \bP(a<Q_n<b) &= (b-a)10^kN_k\qquad \left(0\ls a < b \ls \frac{1}{10^kN_k} \right) \quad \text{ and } \\  \bP\Big(r_n = \frac{i}{10^kN_k}\Big) &= \frac{1}{10^kN_k} \qquad  \left(0 \ls i < 10^kN_k\right).\end{split} \]
	In other words, the $Q_n$'s follow the continuous uniform distribution in the interval $[0,1/(10^kN_k))$ while the $r_n$'s are discretely uniformly distributed among the points $i/(10^k N_k)$. If we set $t_n = Q_n + r_n \, (n\gs 1)$ then $(t_n)_{n\in \bN}$
	is a sequence of independent random variables that are  uniformly distributed on $[0,1)$. \par The first statement we prove has to do with the gap distribution of $\rn$. 
	\begin{lemma} \label{rn}
		The sequence $\rn$ of random variables defined above has exponential gaps almost surely.
	\end{lemma}
	\begin{proof} Fix some real number $s>0.$ We want to find the asymptotic behavior of
		\begin{equation*}  
			P_N(s)  = \frac{1}{N}\#\Big\{ i\ls N : r_{i+1,N} - r_{i,N} \ls \frac{s}{N} \Big\}.  
		\end{equation*}  
		Given $N \gs 1$ let $k\gs 1$ be such that $N_k < N \ls N_{k+1}$. 
		Consider the following sets of indices $i\ls N$: 
		\begin{align*} 
			A_N &= \left\{ i \ls N : r_{i,N} \text{ or } r_{i+1,N} \text{ comes from one of  } ( r_n)_{n =1}^{2N_{k-1}}  \right\}, \\  
			B_N &=  \left\{ i \ls N : r_{i,N} \text{ and } r_{i+1,N} \text{ both come from } (r_n)_{n = 2N_{k}+1}^{N}  \right\}, \\
			\Gamma_N & =  \left\{ i \ls N, i \notin A_N : \text{ either } r_{i,N} \text{ or } r_{i+1,N} \text{ comes from } (r_n)_{n = 2N_{k-1}+1}^{2N_k}  \right\}, \\ 
			\Delta_N & = \left\{ i \ls N : r_{i,N} \text{ and } r_{i+1,N} \text{ both come from } (r_n)_{n = 2N_{k-1}+1}^{2N_k}  \right\}.
		\end{align*} These sets clearly form a partition of the set of indices $i\ls N.$ 
		
		In order to proceed with the proof, we introduce some more notation. Given an interval $J\subseteq \bN$ of positive integers, we write \begin{equation} \label{xiJ}
			x_{i,J} \qquad (1\ls i \ls |J| )     
		\end{equation}
		for the points $x_n, \,\, n\in J$, written in increasing order; in other words, we have $\{x_{i,J} : 1\ls i \ls |J| \} = \{x_n : n\in J \}$ and $x_{1,J} \ls \ldots \ls x_{|J|, J}.$  We will apply the notation for the ranges 
		\[ F_N = (2N_{k-1}, N]  \qquad \text{ and } \qquad  J_N = (2N_k, N]  \]
		of the index set.
		
		We shall distinguish two cases regarding the position of the value of $N$ inside the range $N_k < N \ls N_{k+1}$. This position will determine which of the terms $A_N,B_N,\Gamma_N,\Delta_N$  give some non-negligible contribution to $P_N(s).$ 
		In both cases, we will prove that the discretization grid is fine enough such that for most points $n \ls N$, we have $Q_n = o(1/N)$ and thus, the almost sure exponential gap distribution of $(t_n)_{n \in \N}$ carries over to the gap distribution of $(r_n)_{n \in \mathbb{N}}$. \newline 
		\noindent Case I. If $ N > 10^{k/2} N_k,$ then 
		\[\begin{split} \frac{|A_N|}{N} &\ls \frac{4N_{k-1}}{N}  < \frac{4}{10^{k/2}},  \quad \frac{|\Gamma_N|}{N} \ls \frac{4N_k}{N}  < \frac{4}{10^{k/2}},\\  \frac{|\Delta_N|}{N} &\ls \frac{2N_k}{N}  < \frac{2}{10^{k/2}}, \quad \frac{|B_N|}{N} = \frac{\lvert J_N\rvert}{N} \gs 1 - \frac{10}{10^{k/2}}\cdot \end{split}\]  
		Thus it suffices to prove that
		\begin{align}\label{pjn}P_{J_N}(s) &:= \frac{1}{ \lvert J_N \rvert} \#\left\{i \ls \lvert J_N \rvert: r_{i+1,J_N} - r_{i,J_N} \ls \frac{s}{\lvert J_N \rvert}\right\} \\ &= \int_0^{s+\varepsilon} e^{-\tau} \mathrm{d}\tau + o(1), \qquad N \to \infty. \nonumber \end{align}
		Define the intervals 
		\[ I_j = \left[\frac{j}{10^{k+1} N_{k+1}},\frac{j+1}{10^{k+1}N_{k+1}}\right), \qquad 0\ls j < 10^{k+1}N_{k+1}\] 
		and given $s>0$ and $0 \ls \ell, j < 10^{k+1} N_{k+1}$  set
		\[ \begin{split}G_N(\ell,j,s) = \#\Big\{i\ls |J_N| : \,\, t_{i+1, J_N} - t_{i, J_N} \ls s/\lvert J_N\rvert\,\,\, \& \,\,\, t_{i+1, J_N},  t_{i, J_N} \in I_j \cup I_{\ell}
			\Big\}
		\end{split}\]
		and 
		\begin{equation}\tilde{G}_N(\ell,j,s) = \#\Big\{ i\ls |J_N| : \,\, r_{i+1, J_N }- r_{i, J_N} \ls s/\lvert J_N\rvert\,\,\, \& \,\,\,
			r_{i+1, J_N}, r_{i, J_N} \in I_j \cup I_{\ell} \Big\}.  \end{equation}
		Now the quantity $P_{J_N}$ defined in \eqref{pjn}	satisfies  \begin{equation}\label{rewrite_pjn}
			P_{J_N}(s) = \frac{1}{\lvert J_N \rvert} \sum_{0\ls j \ls \ell < 10^{k+1} N_{k+1}} \hspace{-3mm}\tilde{G}_N(\ell,j,s).  
		\end{equation}
		{\it Claim:} For any $\varepsilon > 0$ and for all $N\gs 1$ sufficiently large, we have
		\begin{equation}\label{tilde_doesnt_matter_}G_N(\ell,j,s-\varepsilon) \ls \tilde{G}_N(\ell,j,s) \ls G_N(\ell,j,s+\varepsilon).\end{equation}
		{\it Proof of Claim:} Let $\{t_{n_1},t_{n_2},\ldots, t_{n_{t_j}}\},$ $\{t_{m_1},t_{m_2},\ldots, t_{m_{t_\ell}}\}$ denote those elements among the points $t_n,\,\, n \in J_N$ that lie in $I_{j}$ and $I_{\ell},$ respectively. By definition of $\rn$ and $(t_n)_{n \in \N}$, among the points  $r_n \in J_N$ the ones that lie in $I_j$ and $I_{\ell}$ are precisely
		$\{r_{n_1},r_{n_2},\ldots, r_{n_{t_j}}\}$ and $\{r_{m_1},r_{m_2},\ldots, r_{m_{t_\ell}}\},$ respectively.\par
		We now want to compare the gap length $s/\lvert J_N \rvert$  with the distances between endpoints of $I_j$ and $I_\ell.$ We assume without loss of generality that $\ell \gs j.$ \par If $s/\lvert J_N \rvert > (\ell - j +1)/(10^k N_k),$ then \eqref{tilde_doesnt_matter_} is true because
		\[ \tilde{G}_N(\ell,j,s) = G_N(\ell,j,s) =\begin{cases} t_j + t_\ell - 1, &\text{ if } j \neq \ell, \\   t_j - 1, &\text{ if } j = \ell. \end{cases} \]  Similarly, if $s/\lvert J_N\rvert < ( \ell - j)/(10^kN_k)$, then $\tilde{G}_N(\ell,j,s) = G_N(\ell,j,s) = 0$, so the remaining case is when $ (\ell - j )/(10^kN_k) < s/\lvert J_N \rvert < ( \ell - j +1)/(10^k N_k)$. Observe that for $n \in J_N$  we have 
		\begin{equation}\label{grid_level} 0 \ls t_n - r_n = Q_n \ls \frac{1}{10^{k+1}N_{k+1}} = o\Big( \frac{1}{\lvert J_N \rvert} \Big), \qquad N\to \infty \end{equation}
		(since $\lvert J_N \rvert \ls N \ls N_{k+1}$) and we infer that 
		\[ 0 = G_N\Big(\ell,j,s - \frac{\lvert J_N \rvert}{10^{k+1}N_{k+1}}\Big) \ls \tilde{G}_N(\ell,j,s) \ls t_j + t_{\ell} - 1 =
		G_N\Big(\ell,j,s + \frac{\lvert J_N \rvert}{10^{k+1}N_{k+1}}\Big).
		\]
		Since $\lvert J_N \rvert/(10^{k+1} N_{k+1}) = o(1),\,\, N\to\infty,$ this shows \eqref{tilde_doesnt_matter_} and the claim is proved. 
		
		\par We can now combine \eqref{rewrite_pjn} with \eqref{tilde_doesnt_matter_} to deduce that 
		\begin{align*}
			P_{J_N}(s) &\ls  \frac{1}{\lvert J_N \rvert}\sum_{0\ls j \ls \ell < 10^{k+1} N_{k+1}} \hspace{-3mm}G_N(\ell,j,s+\varepsilon) + o(1) \\ &= \frac{1}{|J_N|}\#\Big\{n \ls |J_N| : t_{n+1, J_N } - t_{n, J_N } \ls \frac{s+\varepsilon}{\lvert J_N\rvert} \Big\} +o(1), \qquad N\to \infty  \end{align*} and since $(t_n)_{n \in \N}$ has almost surely exponential gaps, this implies that almost surely,
		\[ \limsup_{N\to \infty}P_{J_N}(s) \ls \int_0^{s+\varepsilon}e^{-\tau}\mathrm{d}\tau. \] A similar argument gives the necessary lower bound for  $\liminf\limits_{N\to \infty}P_{J_N}(s)$. Since $\varepsilon$ can be arbitrarily small, this concludes the proof for the first case.\\
		
		\noindent Case II. If $N\ls 10^{k/2}N_k$, then only $A_N$ is negligible and thus we need to show \eqref{pjn} with $J_N$ replaced by $F_N$.
		Writing now \[
		\begin{split}
			G_N(\ell,j,s) = \#\Big\{i\ls |F_N| : \,\, t_{i+1, J_N} - t_{i, F_N} \ls s/\lvert F_N\rvert\,\,\,  \& \,\,\, t_{i+1, F_N},  t_{i, F_N} \in I_j \cup I_{\ell}
			\Big\}, \\
			\tilde{G}_N(\ell,j,s) = \#\Big\{ i\ls |F_N| : \,\, r_{i+1, F_N }- r_{i, F_N} \ls s/\lvert F_N\rvert\,\,\, \& \,\,\,
			r_{i+1, F_N}, r_{i, F_N} \in I_j \cup I_{\ell}
			\Big\},\end{split}\]
		we can prove once more \eqref{tilde_doesnt_matter_} by replacing \eqref{grid_level} with
		\[
		0 \ls t_n - r_n = Q_n \ls \frac{1}{10^{k}N_{k}} = o\Big( \frac{1}{\lvert F_N \rvert} \Big), \qquad N\to \infty
		\]
		where we used that in this case, $\lvert F_N \rvert \ls N \ls 10^{k/2}N_k$. The rest of the proof follows along the lines of the first case. \end{proof}
	
	Having established that the random sequence $\rn$ has exponential gaps almost surely, we shall now show that the same holds for the random sequence $\xn.$ \par  Given $N\gs 1$ and arbitrary points $w_1,\ldots,w_N$, we shall write
	$$\cG (w_1,\ldots,w_N) = \{w_{i+1,N} - w_{i,N} :  i = 1,\ldots, N\} 
	$$ for the set of gaps among the points $w_1,\ldots, w_N,$ counted with multiplicity. Observe that for the random sequences defined previously, such as $\xn,$  the corresponding sets of gaps (e.g. $\cG(x_1,\ldots, x_N)$) will likewise be random sets.  \par For any $N\gs 1,$ we let $k\gs 1$ be the unique index such that $N_k < N \ls N_{k+1}$. First,  note that since the function $f^{(k)}$ exchanges intervals of the form $[x_{i,N_k}, x_{i+1,N_k})$ we have   
	\begin{equation} \label{set1} \cG(x_1,\ldots, x_N) =  \cG( f^{(k)}(x_1),\ldots, f^{(k)}(x_N)). \end{equation}
	Because by definition $x_n = f^{(k)}(\tilde{y}_n)$ for all $N_k < n \ls 2N_k,$ for these values of $n$ we have $f^{(k)}(x_n) = (f^{(k)}\circ f^{(k)})(\tilde{y}_n) = \tilde{y}_n.$ Thus the set in \eqref{set1} is equal to 
	\begin{equation}\label{set2}
		\cG( f^{(k)}(x_1),\ldots, f^{(k)}(x_{N_k}), \tilde{y}_{N_k+1},\ldots, \tilde{y}_{2N_k}, f^{(k)}(x_{2N_k+1}),\ldots, f^{(k)}(x_N)).  
	\end{equation}
	Observe that here we have assumed that $N>2N_k.$ We will do so in the remaining part of the proof that $\xn$ has exponential gaps almost surely; when $N_k < N \ls 2N_k$ the arguments are similar but much simpler. \par Since $\{ f^{(k)}(x_n) : n\ls N_k\} = \{ x_n : n\ls N_k\}$ the set in \eqref{set2} is equal to 
	\begin{align}
		\cG( x_1,\ldots,x_{N_k}, \tilde{y}_{N_k+1},&\ldots, \tilde{y}_{2N_k}, f^{(k)}(x_{2N_k+1}),\ldots, f^{(k)}(x_N)) = \nonumber \\ =&  \cG( x_1,\ldots,x_{2N_{k-1}}, \tilde{y}_{2N_{k-1}+1},\ldots, \tilde{y}_{2N_k}, f^{(k)}(x_{2N_{k}+1}), \ldots, f^{(k)}(x_N)) \label{set3}  \\
		=& \cG( x_1,\ldots,x_{2N_{k-1}}, \tilde{y}_{2N_{k-1}+1},\ldots, \tilde{y}_{2N_k}, f^{(k)}(\tilde{y}_{2N_{k}+1}), \ldots, f^{(k)}(\tilde{y}_N)). \nonumber
	\end{align} 
	Here we simply used the fact that for $2N_{k-1} < N \ls N_{k}$ and $2N_k < N \ls N_{k+1}$ we have $x_n = \tilde{y}_n.$ If we iterate the same arguments as above with $f^{(k-1)}$ instead of $f^{(k)},$ we see that the set in the right hand side of \eqref{set3} is equal to 
	\begin{multline*}
		\cG( f^{(k-1)}(x_1),\ldots,f^{(k-1)}(x_{2N_{k-1}}), f^{(k-1)}(\tilde{y}_{2N_{k-1}+1}),\ldots, f^{(k-1)}(\tilde{y}_{2N_k}), 
		\\  \qquad \qquad  (f^{(k-1)}\circ f^{(k)})(\tilde{y}_{N_{k}+1}), \ldots, (f^{(k-1)} \circ f^{(k)})(\tilde{y}_N))
		\\ \hfill  = \cG( f^{(k-1)}(x_1),\ldots,f^{(k-1)}(x_{N_{k-1}}), \tilde{y}_{N_{k-1}+1},\ldots, 
		\tilde{y}_{2N_{k-1}},f^{(k-1)}(\tilde{y}_{2N_{k-1}+1}), \ldots)
		\\ \hfill  = \cG( x_1,\ldots,x_{N_{k-1}}, \tilde{y}_{N_{k-1}+1},\ldots, 
		\tilde{y}_{2N_{k-1}},f^{(k-1)}(\tilde{y}_{2N_{k-1}+1}), \ldots)
		\\ \hfill = \cG(x_1,\ldots,x_{2N_{k-2}}, \tilde{y}_{2N_{k-2}+1},\ldots, \tilde{y}_{2N_{k-1}},f^{(k-1)}(\tilde{y}_{2N_{k-1}+1}), \ldots).
	\end{multline*}
	The same procedure can be applied for $f^{(\ell)}, \ell  = k-2,k-3,\ldots,1$ and thus we finally deduce that
	\begin{equation} \label{same_gaps}
		\cG(x_1,\ldots,x_N) = \cG(z_1,\ldots,z_N),
	\end{equation}
	where we set $z_n = \tilde{y}_n$ for $1\ls  n \ls 2N_1$ and
	\[
	z_n = (f^{(1)}\circ \ldots \circ f^{(\ell)}) (\tilde{y}_n) \qquad \text{ if } \quad  2N_{\ell} < n \ls 2N_{\ell +1}.
	\]
	
	In view of \eqref{same_gaps}, it suffices to prove that the sequence $\zn$ has exponential gaps almost surely.
	
	\begin{lemma} \label{same_distribution}
		For any $N\gs 1,$  we have $(z_1,\ldots,z_N)\stackrel{d}{=} (\tilde{y}_1,\ldots, \tilde{y}_N).$
	\end{lemma}
	\begin{proof} For any $k\gs 1$, consider the random vector $ \widetilde{Y}_k = (\tilde{y}_1,\ldots,\tilde{y}_{2N_k}). $   The functions $f^{(1)}, \ldots, f^{(k)}$ defined earlier are random maps, that depend only on the vector $\widetilde{Y}_k.$ Thus given $\ba = (a_1, \ldots, a_{2N_k}) \in \mathbb{R}^{2N_k}$ with all of its coordinates lying inside the set $P_k = \{ i/10^{k+1}N_{k+1} : 0\ls i < 10^{k+1}N_{k+1} \}$,  we can write $f_{\ba}$  for the map $f^{(k)}\circ \cdots \circ f^{(1)}$ when $\widetilde{Y}_k=\ba.$
		Furthermore, note that since each map $f^{(k)}$ is a bijection when restricted to $P_k$ and the random variables  $\tilde{y}_n$ assume each of the values $i/10^{k+1}N_{k+1}$ with equal probability, for all $n> 2N_k$ we have 
		\begin{equation}\label{f_invariant} \bP\left( \tilde{y}_n \in f_{\ba}(A) \right)  = \bP\left( \tilde{y}_n \in A \right),  \quad \text{ for any }  A\subseteq [0,1]. \end{equation}    
		{\it Claim 1:} We have $z_n \stackrel{d}{=} \tilde{y}_n \text{ for all } n \ls N$.  \newline
		{\it Proof of Claim 1.} Let $n\ls N$ and $\ell \gs 1$ be such that $2N_{\ell} < n \ls 2N_{\ell +1}.$ Then by definition, $z_n = (f^{(1)} \circ \cdots \circ f^{(\ell)})(\tilde{y}_n)$ and 
		\[\begin{split}
			\mathbb{P} ( z_n < x ) &= \mathbb{P} \left( (f^{(1)} \circ \cdots \circ f^{(\ell)})(\tilde{y}_n)  < x \right) = \mathbb{P} \left( \tilde{y}_n \in  (f^{(\ell)}\circ \cdots \circ f^{(1)})([0,x))  \right)  \\
			&= \sum_{{\ba}\in P_\ell} \mathbb{P} \left( \tilde{y}_n \in (f^{(\ell)}\circ \cdots \circ f^{(1)})([0,x) )  \,  \land \, \widetilde{Y}_\ell = \ba  \right)
			\\&= \sum_{\ba\in P_\ell } \mathbb{P}  \left(\tilde{y}_n \in  f_{\ba} ([0,x))) \land  \widetilde{Y}_\ell = \ba  \right) = \sum_{\ba \in P_\ell} \mathbb{P}  \left( \tilde{y}_n \in  f_{\ba} ([0,x))   \right) \cdot \mathbb{P}  ( \widetilde{Y}_\ell = \ba ) \\
			&= \sum_{\ba \in P_\ell} \mathbb{P}  (\tilde{y}_n \in   [0,x) ) \cdot \mathbb{P}  ( \widetilde{Y}_\ell = \ba ) = \bP( \tilde{y}_n < x).
		\end{split} \]
		
		\noindent {\it Claim 2: } For any $N\gs 1,$ the random variables $z_1,\ldots, z_N$ are mutually independent. \newline 
		{\it Proof of Claim 2.} 
		We define
		\[\begin{split} &Z_1 = (z_1,\ldots,z_{2N_1}), \quad Z_{\ell+1} = (z_{2N_\ell+1},\ldots,z_{2N_{\ell+1}} ) \qquad (1 \ls \ell < k), \\  &\text{ and  } \quad Z_{k+1} = (z_{2N_k+1},\ldots,z_{N}).
		\end{split} \]
		By the definition of the $z_i$'s, all coordinates of a fixed $Z_{\ell}$ are mutually independent, so we are left to show that $Z_1,\ldots, Z_{k+1}$ are mutually independent. We will use the fact that if $X_1,\ldots, X_{k+1}$ are random variables and for $1 \ls \ell \ls k$, $X_{\ell+1}$ is independent of $(X_1,\ldots, X_{\ell})$, then $(X_1,\ldots, X_{k+1})$ are mutually independent - this can be deduced directly from the definitions. Thus, it suffices to show that  $Z_{\ell+1}$ is independent of $(Z_1,\ldots, Z_{\ell})$ for all $1 \ls \ell \ls k$. We will prove the case $\ell < k$; the case $\ell = k$ works completely analogously. 
		Writing $Y_1 = (\tilde{y}_1,\ldots,\tilde{y}_{2N_1}), \quad Y_{j+1} = (\tilde{y}_{2N_{j}+1},\ldots,\tilde{y}_{2N_{j+1}} ), 1 \ls j \ls \ell,$
		we observe that for any admissible $\ba = (a_1,\ldots,a_{2N_{\ell}}),$ $\bb = (b_{2N_{\ell}+1},\ldots,b_{2N_{\ell+1}})$, we have 
		\begin{align*}
			\bP \left(Z_{\ell+1}=\bb \land (Z_1,\ldots,Z_{\ell}) = \ba \right) &=  
			\bP \left(f_{\ba}^{-1}(Y_{\ell+1}) = \bb \land  (Z_1,\ldots,Z_{\ell}) = \ba \right) 
			\\& = \bP \left(Y_{\ell+1} = f_{\ba}(\bb) \land  (Z_1,\ldots,Z_{\ell}) = \ba \right) ,
		\end{align*}
		where $f_{\ba} := f^{(\ell)}\circ \cdots \circ f^{(1)}$ when $(Z_1,\ldots,Z_{\ell}) = \ba$.
		Since $f_{\ba}(\bb)$ is now a deterministic value and $(Z_1,\ldots,Z_{\ell})$ only depends on $(Y_1,\ldots,Y_{\ell})$, which is independent of $Y_{\ell+1}$, we have
		\[
		\begin{split}
			\bP \left(Y_{\ell+1} = f_{\ba}(\bb) \land  (Z_1,\ldots,Z_{\ell}) = \ba \right)
			&= \bP(Y_{\ell+1} = f_{\ba}(\bb))\cdot \bP ((Z_1,\ldots,Z_{\ell}) = \ba)
			\\&=  \bP(Y_{\ell+1} = \bb)\cdot \bP ((Z_1,\ldots,Z_{\ell}) = \ba)
			\\&= \bP(Z_{\ell+1} = \bb)\cdot \bP ((Z_1,\ldots,Z_{\ell}) = \ba),
		\end{split}
		\]
		where we used \eqref{f_invariant} in the penultimate equality.  \par Having established the two claims, it follows directly that 
		\begin{align*}
			\bP(z_1<a_1, \ldots, z_N < a_N) & = \bP(z_1<a_1) \cdots \bP(z_N < a_N) = \bP(\tilde{y}_1<a_1) \cdots \bP(\tilde{y}_N < a_N) \\ & = \bP(\tilde{y}_1<a_1,\ldots, \tilde{y}_N < a_N)
		\end{align*}
		for any $N\gs 1$ and $a_1,\ldots, a_N \in [0,1],$ which proves of the lemma.
	\end{proof}
	By Lemma \ref{same_distribution}, $\zn$ has the same distribution as $(\tilde{y}_n)_{n\in\bN}.$ Since this distribution is the same as that of the sequence $\rn$ in Lemma \ref{rn}, it follows that $\zn$ has exponential gaps almost surely.  Finally \eqref{same_gaps} allows us to deduce that the sequence $\xn$ has exponential gaps almost surely. 
	
	\subsection{Uniform Distribution}  So far we have not specified what exactly the sequence $(N_k)_{k \in \N}$ is; what we know is that as long as it satisfies \eqref{ratio}, the resulting sequence $\xn$ has exponential gaps almost surely. Now we will show that $(N_k)_{k \in \N}$ can be defined in such a way that $\xn$ fails to be uniformly distributed with positive probability. Accordingly, from now on the terms of $(N_k)_{k\in \bN}$ are not given, but we rather intend to show that they can be chosen inductively such that $\xn$ has the desired properties. Whenever we make the assumption that all or some of the terms of $(N_k)_{k\in \bN}$ have been chosen,  this will be stated explicitly.	\par
	
	We start by recalling some of the parameters in the construction of the sequence $\xn$ in the previous section. For any $k\gs 1$ and $0\ls \ell \ls 10^kN_k,$ assuming $N_1, \ldots, N_k$ have been chosen, we defined 
	\begin{align*}
		\mu_k(\ell) & = \# \{  i\in \MlL : \#\{  N_k < n \ls 2N_k : \tilde{y}_n \in [x_{i,N_k}, x_{i+1,N_k})  \} =0  \}, \\
		\nu_k(\ell) & = \#\{  i\in \MlR :   \#\{  N_k < n \ls 2N_k : \tilde{y}_n \in [x_{i,N_k}, x_{i+1,N_k})  \} >0  \}.
	\end{align*} 
	It will be also useful to write
	\[\begin{split} \TlL &= \big\{\tilde{y}_n: N_k < n \ls 2N_k \text{ and }  \tilde{y}_n \in \cup\{ [x_{i,N_k},  x_{i+1,N_k}) : i\in \MlL  \} \big\}, \\  \TlR &= \big\{\tilde{y}_n: N_k < n \ls 2N_k \text{ and }  \tilde{y}_n \in \cup\{ [x_{i,N_k},  x_{i+1,N_k}) : i\in \MlR  \} \big\} \end{split}  \] for the sets of points $\tilde{y}_n$ that lie within some interval of length $\ell/(10^kN_k)$ for some $i< N_k/2$ and $i\gs N_k/2$, respectively. Furthermore we write $m_k = x_{N_k/2, N_k}$ for the $(N_k/2)$--th point among the $x_n,\,\, 1\ls  n\ls N_k$, when they are written in increasing order.

	With these definitions, we observe that 
	\begin{equation} \label{observation}
		\text{ if $\quad \mu_k(\ell) \gs \nu_k(\ell), \quad$ then
			$\quad f^{(k)}\big(\TlR \big) \subseteq [0,m_k).$}   \end{equation} 
	
	In order to complete the proof of Theorem \ref{thm1}, we need a lemma that  shows that with high probability, the assumption that $\mu_k(\ell) \gs \nu_k(\ell)$ is fulfilled for many $k$ and $\ell$ and thus, according to  \eqref{observation}, the number of points that are shifted from $[m_k,1)$  into $[0,m_k)$ is large (unsurprisingly, $m_k$ will turn out to converge to $1/2$ almost surely). This will imply that the ``shifting'' procedure is applied sufficiently often to result in a sequence $\xn$ which is not uniformly distributed.

	\begin{lemma}[Key Lemma]\label{key_lemma_ud_}
		There exist constants $A,B,C> 0$ with $B > A $ such that for any $\varepsilon> 0$ there exists a sequence $(N_k)_{k\in \N}$ (depending on $\varepsilon)$ that satisfies \eqref{ratio}, and for which additionally the following holds with probability at least $1 - \varepsilon$. Writing 
		$$
		J_k(A,B) =  \Big\{\ell \in \mathbb{N}: \frac{1}{B} \ls \frac{\ell}{10^k} < \frac{1}{A}\Big\}, \qquad k \gs 1,
		$$
		for any $k \in \N$ and any $\ell \in J_k(A,B)$, there exists $j(\ell) \in \{2\ell, 2\ell+1\}$ with 
		\begin{enumerate}
			\item[(i)] $\mu_k(j(\ell)) \gs \nu_k(j(\ell))$ and
			\item[(ii)] $\lvert \TjlR \rvert \gs C \dfrac{N_k}{10^k}.$
		\end{enumerate}
	\end{lemma}

	The proof of Lemma \ref{key_lemma_ud_} will be given in Section \ref{sec_key_lemma}. For the moment we take the lemma for granted, and show that it implies that the resulting sequence $\xn$ is not uniformly distributed with positive probability. 
	We claim that for the sequence $(N_k)_{k \in \bN}$ as in Lemma \ref{key_lemma_ud_}, we have $\lim_{k \to \infty} m_k = \frac{1}{2}$  with probability $1$. Indeed, assume there exists some $\delta>0$ such that \begin{equation} \label{large} m_k \ls \frac{1}{2} - 2\delta \qquad \text{ for inf. many } k\gs 1. \end{equation}
	Then for all sufficiently large values of $k$ for which \eqref{large} holds we have $x_{1,N_k}, \ldots, x_{N_k/2, N_k}$ $ \ls \frac12 - 2\delta $ and 
	\begin{align*}
		\#\big\{ n \ls N_k : y_n \ls \tfrac12 - \delta\big\} &\gs \#\big\{ 2N_{k-1} < n \ls N_k : y_n \ls \tfrac12 - \delta \big\} \\
		& \gs \#\big\{ 2N_{k-1} < n \ls N_k : x_n \ls \tfrac12 - 2\delta \big\} > \tfrac12 N_k,
	\end{align*}which happens for infinitely many $k$ with probability $0$ because $\yn$ is u.d. mod $1$ almost surely. We argue similarly when $m_k \gs \frac12 + 2\delta$ for some $\delta>0$, which proves $\lim_{k \to \infty} m_k = \frac{1}{2}$.
	
	To complete the proof that $\xn$ is not uniformly distributed, it is sufficient to show that with positive probability, there exists $\delta > 0$ such that 
	\begin{equation}\label{bias1/2}
		\liminf_{k \to \infty} \frac{1}{2N_k}\#\left\{1 \ls n \ls 2N_k: x_n \in [0,m_k)\right\} > \frac{1}{2} + \delta.
	\end{equation}  
	We have 
	\begin{equation}\label{calculation}
		\begin{split}
			\#\left\{ n \ls 2N_k: x_n \in [0,m_k)\right\} 
			= &\#\left\{ n \ls N_k: x_n \in [0,m_k)\right\} 
			\\&+ \#\left\{N_k < n \ls 2N_k: x_n \in [0,m_k)\right\}.
		\end{split}
	\end{equation}
	For the first term appearing on the right hand side of \eqref{calculation} we have
	\begin{equation}\label{calculation2} \#\left\{  n \ls N_k: x_n \in [0,m_k)\right\} =\frac{N_k}{2} - 1 \end{equation} by the definition of $m_k.$ For the second term on the right hand side of \eqref{calculation}, by the definition of the sequence $\xn$ we have 
	\begin{equation}\label{calculation3}
		\begin{split}
			\#\left\{N_k < n \ls 2N_k: x_n \in [0,m_k)\right\} 
			= &\#\left\{N_k < n \ls 2N_k: \tilde{y}_n \in [0,m_k)\right\} 
			\\&+ \#\left\{N_k < n \ls 2N_k: \tilde{y}_n \neq x_n \right\}.
		\end{split}
	\end{equation}

	Since the numbers $(y_n)_{n = N_{k}+1}^{2N_k}$ are independent uniformly distributed random variables, $\tilde{y}_n - y_n \ls 1/(10^k N_k)$, and $m_k \to 1/2$ almost surely, we can deduce that
	\begin{equation}\label{calculation4}\frac{1}{2N_k}\#\left\{N_k < n \ls 2N_k: \tilde{y}_n \in [0,m_k)\right\} 
		= \frac{1}{4} + o(1),\qquad k \to \infty \end{equation} 
	almost surely. Combining \eqref{calculation} with \eqref{calculation2}, \eqref{calculation3}, \eqref{calculation4} we get that 
	\[  \frac{1}{2N_k}\#\left\{ n \ls 2N_k: x_n \in [0,m_k)\right\} \gs \frac{1}{2} + \frac{1}{2N_k}\#\left\{N_k < n \ls 2N_k: \tilde{y}_n \neq x_n \right\}, \quad k\to \infty,   \]
	and it remains to prove that 
	\begin{equation} \label{remains}
		\frac{1}{2N_k} \#\left\{N_k < n \ls 2N_k: \tilde{y}_n \neq x_n \right\} \gs  \delta  \end{equation}
	holds with positive probability for some fixed $\delta > 0$. From \eqref{observation} we  deduce that
	\[\#\left\{N_k < n \ls 2N_k: \tilde{y}_n \neq x_n \right\} \gs
	\sum_{\substack{1\ls \ell \ls 10^k N_k \\ \mu_k(\ell) \gs \nu_k(\ell)}} \lvert \TlR \rvert.
	\]
	We apply Lemma \ref{key_lemma_ud_}: for any $10^k/B \ls \ell \ls 10^k/A$, with probability at least $1 - \varepsilon$, there exists $j(\ell) \in \{2\ell, 2\ell+1\}$ such that $\mu_k(j(\ell)) \gs \nu_k(j(\ell))$ and $\lvert \TjlR\rvert \gs C N_k/10^k.$ Hence with probability at least $1-\varepsilon,$
	\[\begin{split}\#\left\{N_k < n \ls 2N_k: \tilde{y}_n \neq x_n\right\} &\gs \sum_{\substack{1\ls \ell \ls 10^k N_k\\ \mu_k(\ell) \gs \nu_k(\ell)}} \lvert \TlR \rvert  
		\gs \sum_{ \frac{10^k}{B} \ls \ell \ls \frac{10^k}{A}} \lvert \TjlR \rvert
		\\&\gs \sum_{ \frac{10^k}{B} \ls \ell \ls \frac{10^k}{A}} C\frac{N_k}{10^k}  =   \left\lfloor 10^k \frac{B-A}{AB} \right\rfloor C\frac{N_k}{10^k}
		\\&\gs 2\delta N_k
	\end{split}
	\]
	for $\delta > 0$ chosen sufficiently small, thereby establishing \eqref{remains}.

	\section{Proof of Lemma \ref{key_lemma_ud_}}\label{sec_key_lemma}
	
	The proof of Lemma \ref{key_lemma_ud_} is given in two steps. In the first step (Lemma \ref{munu_lem_}) we prove that  when the index $\ell$ satisfies a certain regularity condition (this is condition \eqref{approx_exp} below), then {\rm (i)} and {\rm (ii)} of the conclusion of Lemma \ref{key_lemma_ud_} hold.  In the second step (Lemma \ref{pigeonhole}) we show that with high probability the regularity condition is indeed fulfilled for many values of $\ell$.  \par Throughout this section, we use the following notation. For any two integers $M <  N$  we write $x_{i,(M,N]}$ for the $i$--th point in increasing order among the points $(x_n)_{n=M+1}^{N}$ as in \eqref{xiJ}.  For convenience we also write \[ J_k = (2N_{k-1}, N_k] \qquad \text{and} \qquad   S_k = N_k - 2N_{k-1}\qquad (k\gs 1). \]  
	Since we will focus our attention on the points $x_n$ for $2N_{k-1}< n \ls N_k$, 
	we define the sets
	\[\tMl  = \Big\{ i\ls S_k: x_{i+1,J_k} -  x_{i, J_k} = \frac{\ell}{10^kN_k} \Big\}\]
	and
	\begin{eqnarray*} \tMlL  & = & \Big\{ i \in  \tMl: i \ls S_k/2 \Big\}, \\
		\tMlR & = & \Big\{ i \in  \tMl: i > S_k/2 \Big\}.
	\end{eqnarray*}
	We begin by calculating the probability distribution of the random variable $|\tMlL|$ under the condition that $|\tMl|=r$ for some positive integer $r$. 
	
	\begin{lemma} \label{prob_distribution}
		Let $k \gs 1$ and $N_1, \ldots, N_{k}$ be given. Let $0< r < S_k.$ With the notation as above, for $0\ls j \ls r$ we have 
		\begin{equation} \label{probjj1}
			\prob\Big( |\tMlL| = j \,\,\big\vert\,\, |\tMl| = r  \Big) =  \binom{S_k}{r}^{-1} \binom{S_k/2}{j}\,\,\,\binom{S_k/2} {r-j}. \end{equation}
	\end{lemma}
	\begin{proof} Recall that by definition we have   
		\[x_n = \frac{1}{10^k N_k}\lfloor 10^k N_k y_n \rfloor, \qquad  2N_{k-1} < n \ls N_k, \] where $(y_n)_{2N_{k-1} < n \ls N_k}$  are independent uniformly distributed random variables.   Thus if we fix a set $\{i_1< i_2 < \ldots < i_r\} \subseteq \{1,\ldots,S_k\}$ of distinct integers, then 
		\begin{equation*} \prob\Big( \tMl = \{i_1 < i_2 < \ldots < i_r \}  \Big\vert \,\, |\tMl|=r \Big) = {\frac{1}{\binom{S_k}{r}}} \cdot \end{equation*} 
		It remains to count subsets of $\{1,\ldots,S_k\}$ with exactly $j$ elements in $\{1,\ldots,S_k/2\}$ and $r-j$ elements in $\{S_k/2+1,\ldots,S_k\}$. As this number is  $\binom{S_k/2}{j}\cdot\binom{S_k/2} {r-j}$, we deduce \eqref{probjj1}.   
	\end{proof} 
	
	We now prove  that with high probability, both $\MlL$ and $\MlR$ contain at least a fixed positive proportion of the elements of $\Ml.$  This should not be a big surprise, since under the assumption of ``randomness'' one would expect roughly as many gaps of the same size among the first $N_k/2$ gaps as among the last $N_k/2$ gaps. 
	
	\begin{lemma}\label{left_eq_right_}
		Let $k \gs 1$ and $N_1, \ldots, N_{k-1}$ be given. For any $\varepsilon > 0$, there exist numbers $C_k = C_k (\varepsilon,N_1, \ldots, N_{k-1}) > 0$ and $M_k = M_k(\varepsilon,N_1, \ldots, N_{k-1}) > 0$ such that for any $N_k \gs M_k$, for any  $1 \ls \ell < 10^k$, and for any  $C_k \ls r \ls S_k$ we have
		\begin{equation}\label{no_bias_lr_}\prob\Big( \min \left\{|\MlL|,|\MlR| \right\} < \tfrac{1}{11}|\Ml|\,\, \Big\vert \,\, |\tMl| = r \Big) \ls \varepsilon.\end{equation}
	\end{lemma}
	\begin{proof}
		We will actually show that when $C_k>0$ and $N_k$ are both large enough, we have 
		$$
		\prob\Big( |\tMlL| < \tfrac{1}{10} |\tMl|\, \Big\vert \, |\tMl| = r  \Big) \ls \frac{\varepsilon}{2}
		$$
		and
		$$
		\prob\Big( |\tMlR| < \tfrac{1}{10}|\tMl| \, \Big\vert \, |\tMl| = r  \Big) \ls \frac{\varepsilon}{2}
		$$
		for all $r$ in $C_k \ls r \ls S_k$. This is sufficient to deduce \eqref{no_bias_lr_},  because by definition the quantities
		\[  |\MlL| -|\tMlL|, \quad  |\MlR|-|\tMlR|, \quad |\Ml | - |\tMl| \]
		are all non-negative and bounded above by $2 N_{k-1}$, and so assuming that $C_k$ is chosen so large such that $2 N_{k-1} \ls \frac{C_k}{10}$ we have 
		\begin{multline*}  \prob\Big( \min \left\{ |\MlL|,|\MlR| \right\} < \tfrac{1}{11}|\Ml| \,\, \Big\vert \,\, |\tMl| = r \Big) \ls  \\
			\hfill \ls  \prob\Big( \min \left\{ |\tMlL|,|\tMlR| \right\} < \tfrac{1}{11}|\Ml|\,\, \Big\vert \,\, |\tMl| = r \Big)  \\
			\hfill \ls  \prob\Big( \min \left\{ |\tMlL|,|\tMlR| \right\} < \tfrac{1}{11}(|\tMl|+2N_{k-1})\,\, \Big\vert \,\, |\tMl| = r \Big) \\
			\hfill \ls  \prob\Big( \min \left\{ |\tMlL|,|\tMlR| \right\} < \tfrac{1}{11} \left(|\tMl|+\frac{|\tMl|}{10} \right)\,\, \Big\vert \,\, |\tMl| = r \Big) \\
			\hfill \ls  \prob\Big( |\tMlL| < \tfrac{1}{10} |\tMl|\, \Big\vert \, |\tMl| = r  \Big) + \prob\Big( |\tMlR| < \tfrac{1}{10}|\tMl| \, \Big\vert \, |\tMl| = r  \Big),
		\end{multline*}
		where we used that $r \gs C_k$ by assumption. \par Assume now that $|\tMl|=r$ with  $C_k \ls r \ls S_k$ is fixed. Applying Lemma \ref{prob_distribution}, we see that for $j < (r-1)/2$
		\begin{align*} \frac{ \prob\left( |\tMlL| = j+1 \,\,\Big\vert\,\, |\tMl|= r \right)}{ \prob\left( |\tMlL| = j \,\,\Big\vert\,\, |\tMl|= r \right)}
			&= \frac{(\frac{S_k}{2}-j)\cdot (r-j)}{(\frac{S_k}{2}-r+j+1)\cdot (j+1)} \gs \frac{r-j}{j+1} \cdot 
		\end{align*}
		Iterating the previous inequality leads to 
		\begin{multline*}
			\prob\Big( |\tMlL| = \Big\lceil \frac{r}{10}\Big\rceil \,\,\Big\vert\,\, |\tMl|= r \Big) \ls \\
			\hfill \ls  \Bigg(\dfrac{2\left\lceil \frac{r}{10}\right\rceil +1}{r-2\left\lceil \frac{r}{10}\right\rceil}\Bigg)^{\left\lceil \frac{r}{10}\right\rceil }\cdot \prob\Big(|\tMlL| = 2\Big\lceil \frac{r}{10}\Big\rceil \,\,\Big\vert\,\, |\tMl|= r \Big) \\
			\hfill \ls  {3}^{-\frac{r}{10}} \cdot \prob\Big( |\tMlL| = 2\Big\lceil \frac{r}{10}\Big\rceil \,\,\Big\vert\,\, |\tMl|= r \Big),
		\end{multline*} 
		provided $C_k>0$ is large enough. Since \eqref{probjj1} implies that  for $j < (r-1)/2$ 
		\[\prob\Big( |\tMlL| = j+1  \,\,\Big\vert\,\, |\tMl|= r  \Big) \gs\prob\Big( |\tMlL| = j  \,\,\Big\vert\,\, |\tMl|= r  \Big),\] we have
		\begin{eqnarray*}
			\prob\Big( |\tMlL| < \frac{r}{10} \,\,\Big\vert\,\, |\tMl|= r  \Big)
			& \ls & \frac{r}{10}~ \prob\Big( |\tMlL| = \Big\lceil \frac{r}{10}\Big\rceil \,\,\Big\vert\,\, |\tMl|= r  \Big)
			\\ & \ls &  \frac{r {3}^{-\frac{r}{10}}}{10} ~ \prob\Big( |\tMlL| = 2\Big\lceil \frac{r}{10}\Big\rceil \,\,\Big\vert\,\, |\tMl|= r \Big)
			\\ & \ls & \frac{r {3}^{-\frac{r}{10}}}{10}.
		\end{eqnarray*}
		Since by assumption we only consider $r$ such that $r \gs C_k$, choosing $C_k=C_k(\varepsilon)$ sufficiently large yields
		\[  \prob\Big( |\tMlL| < \frac{r}{10} \,\,\Big\vert\,\, |\tMl|= r \Big) \ls \frac{\varepsilon}{2}.\]
		By symmetry the same holds for $|\tMlR|$ instead of $|\tMlL|$, so the desired result follows.
	\end{proof}
	With Lemma \ref{left_eq_right_} proved, we can proceed to the first step towards the Key Lemma \ref{key_lemma_ud_}. We show that it is very unlikely that conditions {\rm (i)} and {\rm (ii)} in its conclusion are not both satisfied when $\lvert  \Ml \rvert$ is of ``typical'' size (more precisely, when $\ell \gs 1$ fulfills the regularity condition \eqref{approx_exp} in the statement of the lemma below). 
	
	\begin{lemma}\label{munu_lem_}
		Let $k\gs 1$ and $N_1,\ldots,N_{k-1}$ be given. 
		There exists an absolute constant $A' > 0$ such that the following holds: 
		For any $B' > A'$ and any $\varepsilon > 0$ there are positive numbers $C = C(A',B')$ and $M = M(A',B',\varepsilon)$ such that
		for any $N_k \gs M$ and any $\ell \in J_k(A',B')$ the probability that the event
		\begin{equation}\label{approx_exp}\frac{1}{4}e^{-\frac{\ell}{10^k}}\frac{N_k}{10^k} \ls \lvert  \tMl \rvert \ls 4e^{-\frac{\ell}{10^k}}\frac{N_k}{10^k}, \end{equation}
		occurs, but that simultaneously at least one of the two events
		\begin{enumerate}
			\item[(i)] $\mu_k(\ell) \gs \nu_k(\ell)$,
			\item[(ii)] $\lvert \TlR \rvert \gs C \dfrac{N_k}{10^k}$
		\end{enumerate}
		does not occur, is at most $\varepsilon$.
	\end{lemma}
	\begin{proof}
		Formally speaking, we try to bound the probability $ \prob \left( \eqref{approx_exp} \land \lnot ( {\rm (i)} \land {\rm (ii)} \right). $
		We split this in the form
		\begin{multline*}
			\prob \left( \eqref{approx_exp} \land \lnot ( {\rm (i)} \land {\rm (ii)} ) \land \left( \min \big\{|\MlL|,|\MlR| \big\} \gs \tfrac{1}{11} |\tMl| \right) \right) \\
			+ \prob \left( \eqref{approx_exp} \land \lnot ( {\rm (i)} \land {\rm (ii)} ) \land \left( \min \big\{|\MlL|,|\MlR| \big\} < \tfrac{1}{11} |\tMl| \right) \right).
		\end{multline*}
		By Lemma \ref{left_eq_right_}, for sufficiently large $N_k$ the expression in the second line of the equation above is bounded by 
		\begin{multline*}
			\prob \left( \eqref{approx_exp} \land \lnot ( {\rm (i)} \land {\rm (ii)} ) \land \left( \min \big\{|\MlL|,|\MlR| \big\} < \tfrac{1}{11} |\tMl| \right) \right) \\
			\ls  \prob \left( \eqref{approx_exp}  \land \left( \min \big\{|\MlL|,|\MlR| \big\} < \tfrac{1}{11} |\tMl| \right) \right) \\
			= \sum \prob \left( \left(\lvert \tMl \rvert = r \right) \land \left( \min \big\{|\MlL|,|\MlR| \big\} < \tfrac{1}{11} |\tMl| \right) \right) \\
			= \sum \prob \left( \min \big\{|\MlL|,|\MlR| \big\} < \tfrac{1}{11} |\tMl| \,\, \Big\vert \,\, |\tMl| = r \right) \prob \left( \lvert \tMl \rvert = r \right) \\
			\ls  \sum \varepsilon ~\prob \left( \lvert \tMl \rvert = r \right)  \ls \varepsilon,
		\end{multline*}
		where all summations in the equation above are understood to be taken over the variable $r$ in the range $\frac{1}{4}e^{-\frac{\ell}{10^k}}\frac{N_k}{10^k} \ls r \ls 4e^{-\frac{\ell}{10^k}}\frac{N_k}{10^k}$. We used that $\ell/10^k \ls 1/A'$ since we assume that $\ell \in J_k(A',B')$, and thus when we choose $M$ sufficiently large, for any $N_k \gs M$, we have $\frac{1}{4}e^{-\frac{\ell}{10^k}}\frac{N_k}{10^k} \gs C_k$ with $C_k$ as in 
		Lemma \ref{prob_distribution}.\\
		
		It remains to deal with
		$$
		\prob \left( \eqref{approx_exp} \land \lnot ( {\rm (i)} \land {\rm (ii)} ) \land \left( \min \big\{|\MlL|,|\MlR| \big\} \gs \tfrac{1}{11} |\tMl| \right) \right),
		$$
		which is bounded above by
		\begin{align}
			\prob & \left( \eqref{approx_exp}  \land \lnot  {\rm (i)} \land \left( \min \big\{|\MlL|,|\MlR| \big\} \gs \tfrac{1}{11} |\tMl| \right) \right) \nonumber \\ 
			& \quad +  \prob \left( \eqref{approx_exp} \land \lnot {\rm (ii)} \land \left( \min \big\{|\MlL|,|\MlR| \big\} \gs \tfrac{1}{11} |\tMl| \right) \right) \label{first_this}   \\
			\ls &  \prob \left(\lnot  {\rm (i)} \land   |\MlL| \gs  \tfrac{1}{44}  e^{-\frac{\ell}{10^k}}\frac{N_k}{10^k}   \right)  
			+  \prob \left( \lnot {\rm (ii)} \land |\MlR| \gs  \tfrac{1}{44}  e^{-\frac{\ell}{10^k}}\frac{N_k}{10^k}  \right).  \nonumber
		\end{align}
		We will show later that 
		\begin{equation}\label{many_empty}\lvert \Tl(L)\rvert + \lvert \Tl(R) \rvert  \ls \lvert \Ml(L)\rvert  \quad \text{implies} \quad \mu_k(\ell) \gs \nu_k(\ell).\end{equation}
		Assuming this implication for the moment to be true, we see that the first term in the right hand side of \eqref{first_this} is bounded above by
		\begin{eqnarray}
			\prob \left( \lvert \Tl(L)\rvert + \lvert \Tl(R) \rvert  > \tfrac{1}{44}  e^{-\frac{\ell}{10^k}}\frac{N_k}{10^k}   \right) \label{first_this_2}
		\end{eqnarray}
		We define the sets \[ \Gl^{(k)} = \bigcup_{i \in \Ml} [x_{i,N_k}, x_{i+1,N_k}) \qquad (\ell  \gs 1).  \]
		Since for $N_k < n \ls 2N_k$ we have $y_n \in \Gl^{(k)} $ if and only if $\tilde{y}_n \in \Gl^{(k)},$ it follows that
		\[\begin{split}\lvert \TlL \rvert + \lvert \TlR\rvert
			&=  \#\{N_k < n \ls 2N_k: \tilde{y}_n \in \Gl^{(k)}\}
			\\& = \#\{N_k < n \ls 2N_k: y_n \in \Gl^{(k)}\}.
		\end{split}
		\]
		Note that $\Gl^{(k)}$ is a union of disjoint intervals
		of total length
		\[\lambda (\Gl^{(k)} ) = \lvert \Ml \rvert \frac{\ell}{10^k N_k}
		\ls 4e^{-\frac{\ell}{10^k}} \frac{\ell}{10^kN_k} \frac{N_k}{10^k}\ls \frac{4}{A'}\frac{1}{10^k}e^{-\frac{\ell}{10^k}},
		\]
		since $\ell/10^k \ls 1/A'$. Because $(y_n)_{n = N_k+1}^{2N_k}$ are  i.i.d.\ uniformly distributed random variables, we have with probability at least $1 - \varepsilon/2$ (provided $N_k$ is sufficiently large) that 
		\begin{equation}\label{estim_T}
			\begin{split} \lvert \TlL \rvert + \lvert \TlR\rvert &= 
				\#\{N_k < n \ls 2N_k: y_n \in \Gl^{(k)} \} \\
				& \ls  2N_k \lambda\left(\Gl\right) \ls
				\frac{8}{A'}e^{\frac{-2\ell}{10^k}}\frac{N_k}{10^k} \cdot \end{split}\end{equation}
		For $A' \gs 500$, we have $8/A' \ls 1/44$, so the probability in line \eqref{first_this_2} is at most $\varepsilon/2$. To deduce that $\varepsilon/2$ also is an upper bound for the probability in line \eqref{first_this}, it remains to prove the implication \eqref{many_empty}. Clearly,
		\[\begin{split}\nu_k(\ell) &=  \#\big\{  i\in \MlR :   \#\{  N_k < n \ls 2N_k : \tilde{y}_n \in [x_{i,N_k},x_{i+1,N_k})  \} >0  \big\}  \\&   \ls  \#\big\{N_k < n \ls 2N_k: \exists i \in \MlR: \tilde{y}_n \in [x_{i,N_k},x_{i+1,N_k})\big\} \\
			& = \lvert\TlR\rvert,\end{split}\]
		where the first and last equality are just the corresponding definitions.
		By the same argument,
		\[\begin{split}\mu_k(\ell) &= \#\big\{  i\in \MlL : \#\{  N_k < n \ls 2N_k : \tilde{y}_n \in [x_{i,N_k},x_{i+1,N_k})  \} =0  \big\}
			\\&= \lvert \MlL \rvert - \#\big\{  i\in \MlL : \#\{  N_k < n \ls 2N_k : \tilde{y}_n \in [x_{i,N_k},x_{i+1,N_k})  \} >0  \big\}
			\\&\gs \lvert \MlL \rvert -  \#\big\{N_k < n \ls 2N_k: \exists i \in \MlL: \tilde{y}_n \in [x_{i,N_k},x_{i+1,N_k})\big\} 
			\\&= \lvert \MlL \rvert - \lvert \TlL\rvert.
		\end{split}
		\]
		Combining the previous estimates we have established \eqref{many_empty}, which shows that the first term in the right hand side of \eqref{first_this} is indeed bounded above by $\varepsilon/2$. \par 
		To give an upper bound for the second term in the right hand side of \eqref{first_this}, we use an argument similar to the one that produced the estimate \eqref{estim_T}. We have
		$\lvert \TlR\rvert = \#\{N_k < n \ls 2N_k: y_n \in \Gl^{(k)}(R)\},$ where \[\Gl^{(k)}(R) := \bigcup_{i \in \Ml(R)}[x_{i,N_k},x_{i+1,N_k})\]
		is a union of disjoint intervals
		of total length
		$\lambda\big(\Gl^{(k)}(R)\big) = \lvert \Ml(R) \rvert \dfrac{\ell}{10^k N_k}.$ Thus
		$$
		\prob \left( \lnot {\rm (ii)} \land |\MlR| \gs  \tfrac{1}{44}  e^{-\frac{\ell}{10^k}}\frac{N_k}{10^k}  \right) \ls \prob \left( \lnot {\rm (ii)} \land \lambda\big(\Gl^{(k)}(R)\big) \gs \frac{1}{44 B'}e^{-A'}\frac{1}{10^k}	  \right), 
		$$
		where we used that $\dfrac{1}{B'} \ls \dfrac{\ell}{10^k} \ls \dfrac{1}{A'}$. We set $C= e^{-A'}/(50 B')$. Since $(y_n)_{n = N_k+1}^{2N_k}$ arises from an i.i.d.\ sample, it is very unlikely that $\lvert \TlR\rvert$ (the number of random points contained in $\Gl^{(k)}(R)$) is significantly smaller than the total length of $\Gl^{(k)}(R)$. Thus the probability on the right-hand side of the previous displayed equation is less than $\varepsilon/2$, provided $N_k$ is sufficiently large.  This gives the desired upper bound for the second term in the rhs of \eqref{first_this}, and completes the proof of Lemma \ref{munu_lem_}.
	\end{proof}
	We now proceed to the second step of the proof of Lemma \ref{key_lemma_ud_}. We show that condition \eqref{approx_exp} holds with high probability.

	\begin{lemma}\label{pigeonhole}
		Let $k\gs 1$ and $N_1, \ldots, N_{k-1}$ be given.  For $\ell \gs 1$ and $N_k \in \N$, let $\mathcal{C}_{\ell} = \mathcal{C}_{\ell}(k,N_k)$ denote the event
		\[ \mathcal{C}_{\ell}  = \left\{  \frac{1}{5}e^{-\frac{2\ell}{10^k}}\frac{N_k}{10^k} \ls \lvert  \Mtl \rvert \ls 5e^{-\frac{2\ell}{10^k}}\frac{N_k}{10^k} \quad  \text{or} \quad \frac{1}{5}e^{-\frac{2\ell}{10^k}}\frac{N_k}{10^k} \ls \lvert  \Mtlo \rvert \ls 5e^{-\frac{2\ell}{10^k}}\frac{N_k}{10^k} \right\}.\] 
		Then for any $\varepsilon > 0$, there exists an integer $M = M(\varepsilon,k)$ such that for every $N_k \gs M$, we have
		\[ \prob\left(\bigcap_{\ell =1}^{10^k} \mathcal{C_{\ell}}\right) \gs 1- \varepsilon.\]
	\end{lemma}
	\begin{proof}
		We will show that $\mathbb{P}\Big(\bigcap\limits_{\ell =1}^{10^k} \tilde{\mathcal{C}}_{\ell}\Big) > 1 - \varepsilon,$ where we define
		\[  \tilde{\mathcal{C}}_{\ell}  = \left\{  \frac{1}{4}e^{- \frac{2\ell}{10^k}}\frac{N_k}{10^k} \ls \lvert  \tMtl \rvert \ls 4e^{- \frac{2\ell}{10^k}}\frac{N_k}{10^k} \quad 
		\text{or} \quad \frac{1}{4}e^{- \frac{4\ell}{10^k}}\frac{N_k}{10^k} \ls \lvert  \tMtlo \rvert \ls 4e^{- \frac{2\ell}{10^k}}\frac{N_k}{10^k}   \right\}. \] 
		This suffices, since by \eqref{ratio} we have  $|\Ml| - |\tMl|   \ls 2N_{k-1}$, which can be made arbitrarily small in comparison with $N_k$ by choosing the lower bound $M$ in the statement of the lemma sufficiently large. Now let $1 \ls \ell \ls 10^k$ be fixed. We first observe that whenever $y_{i+1,J_k} - y_{i,J_k} \in \left(\frac{\ell}{10^k N_k},\frac{\ell+1}{10^k N_k}\right)$ then the difference  $x_{i+1,J_k} - x_{i,J_k} $ is either $ \ell /(10^k N_k)$ or $(\ell+1)/(10^k N_k);$ in other words, 
		\begin{equation}\label{gap_lengthsxy} y_{i+1,J_k} - y_{i,J_k} \in \left(\frac{2\ell}{10^k N_k},\frac{2\ell+1}{10^k N_k}\right) \,\, \Rightarrow \,\, i \in \tMtl \cup \tMtlo.\end{equation}
		For any $\eta > 0$ there exists $M$ sufficiently large, such that for $N_k \gs M$ we have $1 - \eta \ls S_k / N_k \ls 1$  and thus  
		\begin{equation}\label{S_kN_k}\begin{split}\frac{1}{S_k} \#&\Big\{ i \ls S_k: y_{i+1, J_k} - y_{i,J_k} \in \Big(\frac{2\ell}{10^k N_k},\frac{2\ell+1}{10^k N_k}\Big) 
				\Big\} \gs \\ \gs &\frac{1}{S_k} \# \Big\{  i \ls S_k:
				y_{i+1,J_k} - y_{i,J_k} \in \Big(\frac{2\ell}{10^k S_k} , \frac{2\ell+1}{10^k S_k}(1-\eta)\Big) \Big\}.\end{split}\end{equation}
		We now use the  fact  that, as noted in the introduction, a sequence of i.i.d.\ uniformly distributed random variables has exponential gap distribution almost surely.  Thus, for any $\delta,\varepsilon > 0$ and any $\ell$ in the range $1 \ls \ell \ls 10^k$, there exists some $M = M(\varepsilon,\delta,\ell) \gs 1$ such that for all $N_k \gs M$ (and thus, $S_k$ sufficiently large) we have
		\[ \prob \Bigg(\frac{1}{S_k} \# \Big\{  i \ls S_k:
		y_{i+1,J_k} - y_{i,J_k} \in \Big(\frac{2\ell}{10^k N_k},\frac{2\ell+1}{10^k N_k}\Big) 
		\Big\} \gs \hspace{-5mm}\int\limits_{\left(\frac{2\ell}{10^k},\frac{2\ell +1}{10^k}(1 - \eta)\right)} \hspace{-8mm} e^{- \tau} \,\mathrm{d}\tau - \delta \Bigg) > 1- \frac{\varepsilon}{10^k}. \]
		Choosing $\eta,\delta>0$ sufficiently small, there exists $M = M(\varepsilon)$ such that for $N_k \gs M$ we have
		\[ \prob \Big( \# \Big\{  i \ls S_k:
		y_{i+1,J_k} - y_{i,J_k} \in \Big(\frac{2\ell}{10^k N_k},\frac{2\ell+1}{10^k N_k}\Big) 
		\Big\} \gs \frac{3S_k}{4\cdot 10^k} e^{-\frac{2\ell}{10^k}}\quad  \forall  \ell \ls 10^k \Big) > 1- \varepsilon. \]
		Since for large $M$ we have $3S_k/4 > N_k/2,$  applying \eqref{gap_lengthsxy} we get
		\[\prob \Big( \lvert \tMtl \cup \tMtlo\rvert \gs  \frac{N_k}{2}\frac{1}{10^k}e^{-\frac{2\ell}{10^k}} \quad \forall  \ell \ls 10^k \Big) > 1-\varepsilon.\]
		Clearly, if $\lvert \tMtl \cup \tMtlo\rvert \gs \dfrac{N_k}{2\cdot 10^k} e^{-\frac{2\ell}{10^k}}$, then at least one of the inequalities
		\begin{equation}\label{lower_bound_lem}\lvert \tMtl\rvert \gs \frac{1}{4} \cdot \frac{N_k}{10^k}e^{-\frac{2\ell}{10^k}}, \qquad  \lvert \tMtlo\rvert \gs \frac{1}{4} \cdot \frac{N_k}{10^k}e^{-\frac{2\ell}{10^k}}\end{equation}
		is fulfilled.
		For the upper bound, we argue similarly, only replacing 
		\eqref{gap_lengthsxy} by 
		\[i \in \tMl \quad\Rightarrow\quad y_{i+1,J_k} - y_{i,J_k} \in \left(\frac{\ell-1}{10^k N_k},\frac{\ell+1}{10^k N_k}\right),\]
		which leads for $N_k$ sufficiently large to
		\[ \begin{split}
			\prob \Big(\lvert \tMtl \rvert \ls
			3 e^{-\frac{2\ell}{10^k}}\frac{N_k}{10^k} \text{ and } \lvert \tMtlo \rvert \ls
			3 e^{-\frac{2\ell}{10^k}}\frac{N_k}{10^k}  \quad  \forall \ell \ls 10^k
			\Big) > 1- \varepsilon.
		\end{split} \]
		Combining the previous estimate with \eqref{lower_bound_lem} yields the desired result.
	\end{proof}
	
	We are finally in position to finish the proof of Lemma \ref{key_lemma_ud_}. Let $\varepsilon > 0$ be fixed. We define the terms of the sequence $(N_k)_{k \in \N}$ inductively. We first choose an arbitrary even value for $N_1$. Assume that  for some $k\gs 1,$ the values of $N_1,\ldots, N_{k-1}$ have been chosen. \par  By Lemma \ref{pigeonhole}, for any $\varepsilon_k> 0$  there exists $M_1 = M_1(\varepsilon_k,k)$ such that for any choice of $N_k \gs M_1,$  with probability at least $1 - \varepsilon_k$ the following holds: for any $1\ls \ell \ls 10^k$ there exists $j(\ell) \in \{2\ell,2\ell+1\}$ such that
	\begin{equation}\label{size_of_mjl} \frac{1}{5}e^{-\frac{2\ell}{10^k}}\frac{N_k}{10^k} \ls \lvert  \Mtjl \rvert \ls 5e^{-\frac{2\ell}{10^k}}\frac{N_k}{10^k}. 
	\end{equation}
	Let $A = 3A', B = 2B'$ with $B' > A' > 0$ as in Lemma \ref{munu_lem_}. Then $\frac{10^k}{B} \ls \ell \ls \frac{10^k}{A}$ implies that $\frac{10^k}{B'} \ls j(\ell) \ls \frac{10^k}{A'}.$  Hence by Lemma \ref{munu_lem_}, if $N_k \gs M_2 = M_2(k,\varepsilon_k)$ and additionally \eqref{size_of_mjl} holds, then with probability at least $1 - \varepsilon_k$ we have
	\begin{equation}\label{2nd_property}\mu_k(j(\ell)) \gs \nu_k(j(\ell))\qquad \text{ and }\qquad  \lvert \TjlR \rvert \gs C \frac{N_k}{10^k}.\end{equation}
	Thus, for $N_k \gs \max\{M_1,M_2\}$, we have with probability at least $1 - 2\varepsilon_k$ that both 
	\eqref{size_of_mjl} and \eqref{2nd_property} hold for some $j(\ell) \in \{2\ell,2\ell+1\}$ for all $\ell \in J_k(A,B)$. Setting $\varepsilon_k =  \varepsilon/2^{k+1}$, we now choose $N_k$ such that:
	\begin{itemize}
		\item[(i)] $N_k \gs \max\{M_1(k,\varepsilon_k),M_2(k,\varepsilon_k)\}$,
		\item[(ii)] $N_{k-1} \mid N_k$, and
		\item[(iii)] $N_k \gs k \cdot N_{k-1}$.
	\end{itemize}
	The sequence $(N_k)_{k \in \N}$ defined in this way clearly satisfies \eqref{ratio}. Furthermore, the union bound shows that 
	\[\mu_k(j(\ell)) \gs \nu_k(j(\ell)) \qquad \text{ and } \qquad  \lvert \TjlR \rvert \gs C \frac{N_k}{10^k}\]
	holds for all $k \gs 1$ and  for all $ \ell \in J_k(A,B)$ with probability at least 
	$1 - \sum\limits_{k \gs 1}\dfrac{\varepsilon}{2^k} =  1 - \varepsilon.$ This concludes the proof of Lemma \ref{key_lemma_ud_}. We have already shown in Section \ref{section_2} that Theorem \ref{thm1} follows from Lemma \ref{key_lemma_ud_}, so the proof of Theorem \ref{thm1} is now complete.
	
	\section{Proof of Theorem \ref{gaps_sets_det_distr}}
	Let $\xn$ be a sequence in $[0,1]$ whose gaps are $\cG = \{g_{i,N} : i\ls N, N\gs 1\}.$ We write $G_{i,N} = [\xin, \xiin)$ for the intervals defined by the points $x_1,\ldots, x_N$ and $\mathcal{F}_N = \sigma\left(\{G_{i,N}: i = 1,\ldots,N\}\right)$ for the $\sigma$-algebra generated by these intervals. Also let $\yn \subseteq [0,1]$ be another sequence whose gaps $h_{i,N}= \yiin - \yin$ are the same as those of $\xn,$ in the sense that $\{g_{i,N} : i \ls N \} = \{ h_{i,N} : i \ls N \}$ for all $N\gs 1.$ We will write $H_{i,N} = [y_{i,N}, y_{i+1,N})$ for the intervals formed by the points $y_1,\ldots, y_N,$ analogously to the intervals $G_{i,N}.$ \par For each gap $g_{i,M}$ of the sequence $\xn,$ we define its {\it descendants} at stage $N\gs M$ to be the set 
	\begin{equation} \label{desc} \desc_{N}({ g_{i,M}};\xn)= \{g_{j,N}: 1 \ls j \ls N \text{ and } G_{j,N} \subseteq G_{i,M}\}. \end{equation}
	In other words, $\desc_{N}({ g_{i,M}};\xn)$ is the set of gaps at stage $N$ that appear as lengths of subintervals of $G_{i,M}.$ \par 
	
	We note that the set of descendants as defined above actually does not really depend on the sequence $\xn$, but only on the set $\cG$ of gaps (since the  gaps are distinct by assumption, so that a new gap can only arise by splitting an existing gap). In view of this observation, we write $\desc_N(g_{i,M};\cG)$ for the set defined in \eqref{desc}. For $M,N\gs 1$ with $N > M$ and $k \gs 1$, we set
	\begin{equation} \label{emn} E_{M,N}^{\,\,(k)}(\cG) = \sum_{i = 1}^M  \frac{1}{g_{i,M}}\left(\frac{\# \desc_N(g_{i,M};\cG)}{N}\right)^k.\end{equation} 
	Recall the hypothesis that the sequence $\xn$ has an asymptotic density function $g: [0,1]\to \mathbb{R}.$ We will now establish the connection between $E_{M,N}^{\,\,(k)}(\cG)$ defined earlier and the moments of the density function $g.$
	
	\begin{lemma} \label{E_moment} 
		If $g:[0,1]\to \mathbb{R}$ is the asymptotic density of $\xn$ and  $E_{M,N}^{\,\,(k)}(\cG)$ is as in $\eqref{emn},$ then  
		\begin{equation}  \label{gk}
			\lim_{M \to \infty}\big( \lim_{N \to \infty}E_{M,N}^{\,(k)}(\cG) \big)  = \int_{0}^1 g^k(x) \, \mathrm{d}x \qquad (k\gs 1).
		\end{equation}
	\end{lemma}
	\begin{proof}
		By the  definition of $g$ we have that \[\lim_{N \to \infty}  \frac{\# \desc_N(g_{i,M};\cG)}{N} = \int_{G_{i,M}} g(t) \,\mathrm{d}t.\]
		This implies that 
		\begin{equation} \label{Econverge}
			\lim_{N \to \infty} E_{M,N}^{\,\,(k)}(\cG)   = \sum_{i = 1}^M \frac{1}{g_{i,M}}\Big(\int_{G_{i,M}} g(t) \,\mathrm{d}t\Big)^k
			= \int_0^1 \mathbb{E}[g\vert \mathcal{F}_M]^k(x)\dd x, \end{equation}
		where $\mathbb{E}[g\vert \mathcal{F}_M]$ denotes the conditional expectation of $g$ with respect to the $\sigma$-algebra $\mathcal{F}_M$ (we refer to \cite[Chapter 5]{ergodic} for definitions of the probabilistic terminology used in this paragraph). Note that $\mathcal{F}_N \subseteq \mathcal{F}_{N+1}$ by construction; that is, the sequence of $\sigma$--algebras $(\mathcal{F}_N)_{N \in \mathbb{N}}$ forms a filtration. Note also that the $\sigma$--algebras $\mathcal{F}_N$ converge to the Borel $\sigma$--algebra $\mathcal{B}([0,1])$, because $\xn$ is dense by assumption. Thus, we can apply the martingale convergence theorem in the form of \cite[Thm 5.5]{ergodic} to obtain
		\[\lim_{M \to \infty} \big( \lim_{N \to \infty}E_{M,N}^{\,\,(k)}(\cG) \big) = \lim_{M \to \infty} \int_0^1 \mathbb{E}[g\vert \mathcal{F}_M]^k(x)\, \dd x = \int_{0}^1 g^k(x) \, \mathrm{d}x \]
		for any $k \gs 1$.  \end{proof} 
	
	Since $E_{M,N}^{\,\,(k)}(\cG)$ only depends on the set $\cG$ of gaps, it may refer to either of the sequences $\xn$ or $\yn.$ By Lemma \ref{E_moment} it is reasonable to assume that the limit in \eqref{Econverge} also equals to the $k$-th moment of $\yn.$ The problem is that it is not a priori known whether a limit density exists for $\yn$, or if it has finite $k$-th moments for all $k\gs 1.$ This will be established in the following two lemmas. 
	
	\begin{lemma} \label{lemma_asy}
		Let $F$ be the asymptotic distribution function of $\yn$. Then $F$ is absolutely continuous. 
	\end{lemma}
	
	We note that the asymptotic distribution function $F$ of $\yn$ need not be unique, but we assumed so for convenience (without this assumption, one would need to work along subsequences, but the argument would remain essentially unchanged).  
	
	\begin{proof}[Proof of Lemma \ref{lemma_asy} ]
		Let $\varepsilon > 0$ fixed. Because the function $\int_0^xg(t) \,\mathrm{d}t$ is absolutely continuous, there exists $\delta = \delta(\varepsilon)>0$ such that for any finite collection of disjoint intervals $(I_k)_{k\ls n}$ with $\sum_{k\ls n}|I_k| < \delta$  we have	$\sum_{i \ls n} \int_{I_k} g(x) \,\mathrm{d}x < \varepsilon$. \par Let $(a_k,b_k), \, k =1,\ldots, n $ be disjoint intervals such that $\sum_{k\ls n} (b_k - a_k) < \delta/2.$ 
		Observe that $\yn$ is dense because it has the same gaps as the dense sequence $\xn$.  Thus, there exists a $K = K(\delta) \gs 1$ and a finite index set $J = J(\delta) \subseteq \{1, \ldots, K\}$ with
		\[\bigcup_{k\ls n}(a_k,b_k)\subseteq \bigcup_{i \in J} H_{i,K}\qquad \text{ and } \qquad \sum_{i \in J} h_{i,K} < \delta.\]
		Let $J'= J'(\delta)$ be such that 
		$\{h_{i,K}\}_{i \in J} = \{g_{i',K}\}_{i' \in J'}$ (this exists since the gap structures coincide). This implies
		\begin{multline*}\sum_{k\ls n} (F(b_k) - F(a_k)) \ls \sum_{i \in J}  \lim_{N \to \infty}\left(\frac{\# \desc_N(h_{i,K};\cG)}{N}\right) \\ = \sum_{i \in J'}  \lim_{N \to \infty}\left(\frac{\# \desc_N(g_{i,K};\cG)}{N}\right)
			= \sum_{i \in J'} \int_{G_{i,K}} g(x) \,\mathrm{d}x < \varepsilon,
		\end{multline*} because $\sum_{i\in J'}|G_{i,K}| = \sum_{i\in J}h_{i,K} < \delta.$ Thus $F$ is absolutely continuous.
	\end{proof}
	
	Since $F$ is absolutely continuous by the previous lemma, the sequence $\yn$ has an asymptotic density function $f \in L^1[0,1]$.
	\begin{lemma}
		For any $k\gs 1,$ the function $f^k$ is integrable.
	\end{lemma}
	\begin{proof}
		The same argument that gave \eqref{Econverge} also gives  \begin{equation}\label{Econverge2} \lim_{N \to \infty} E_{M,N}^{\,\,(k)}(\cG)  = \sum_{i = 1}^M \frac{1}{h_{i,M}}\Big(\int_{H_{i,M}} f(t) \,\mathrm{d}t\Big)^k.
		\end{equation}
		If $f^k$ is not integrable, we define the truncations $f_M$ by
		\[f_M(x) = \begin{cases}
			f(x) &\text{ if } f(x) \ls M,
			\\ 0 &\text{ otherwise}.
		\end{cases}\]
		Since clearly 
		\[\sum_{\ell = 1}^M \frac{1}{h_{\ell,M}}\Big(\int_{H_{\ell,M}} f(t) \,\mathrm{d}t\Big)^k \\ \gs \sum_{\ell = 1}^M \frac{1}{h_{\ell,M}} \Big(\int_{H_{\ell,M}} f_M(t) \,\mathrm{d}t\Big)^k  \]
		and $f_M^k$ is integrable, letting $M\to \infty$ in \eqref{Econverge2} and using \eqref{gk} gives
		\[\int_{0}^{1} g^k(x) \,\mathrm{d}x = \lim_{M \to \infty}\big( \lim_{N \to \infty}  E_{M,N}^{\,\,(k)}(\cG) \big) \gs \lim_{M\to \infty} \int_{0}^{1} f_M^k(x) \,\mathrm{d}x =\infty,
		\]
		a contradiction.
	\end{proof}
	
	Having established that $\int_0^1 f^k(x)\dd x < \infty$ for all $k\gs 1,$  the proof of Theorem \ref{gaps_sets_det_distr} now follows from Lemma \ref{E_moment}. Since $E_{M,N}^{\,\,(k)}(\cG)$ is the same for both sequences $\xn$ and $\yn$, and since it converges to the $k$-th moments of the respective asymptotic densities, these moments are necessarily the same for both sequences.


\begin{thebibliography}{HD82}
		\bibitem{abty}
		C. Aistleitner, S. Baker, N. Technau, N. Yesha, {\it Gap statistics and higher correlations for geometric progressions modulo one}, Math. Ann. 385, no.1-2, 845\,--\,861 (2023).
		
		\bibitem{abr}
		C. Aistleitner, V. Blomer, M. Radziwi{\l}{\l}, {\it Triple correlation and long gaps in the spectrum of flat tori}, J. Eur. Math. Soc., to appear. Preprint available at: arXiv 1809.07881 (2018)

		\bibitem{allp} C. Aistleitner, T. Lachmann, P. Leonetti and P. Minelli, {\it On the number of gaps of sequences with Poissonian pair correlations,} Discrete Math. 344, no. 11, Article ID 112555, 13 p. (2021).
		
		\bibitem{alp}
		C. Aistleitner, T. Lachmann and F. Pausinger, {\it  Pair correlations and equidistribution,} J. Number Theory 182, 206\,--\,220, (2018).
		
		\bibitem{ac} 
		D. Altman, Z. Chase, {\it On the smallest gap in a sequence with Poisson pair correlations,} Pre-print: arXiv:2210.17466 (2022).
		
		\bibitem{berry}
		M. V. Berry, M. Tabor, {\it Level clustering in the regular spectrum.} Proc. Roy. Soc. A. 356 (1977), 375\,--\,394.
		
		\bibitem{cy} 
		S. Chaubey, N. Yesha, {\it The distribution of spacings of real-valued lacunary sequences modulo one,} Mathematika 68, no. 2, 416\,--\,428 (2022). 
		
		\bibitem{ergodic}
		M. Einsiedler, T. Ward,  {\em Ergodic Theory: with a view towards Number Theory,} Springer Verlag, London (2011).
		
		\bibitem{elkies}
		N. D. Elkies, C. T. McMullen, {\it Gaps in $\sqrt{n} \bmod 1$ and ergodic theory}, Duke Math. J. 123, no. 1, 95\,--\,139 (2004)
		
		\bibitem{sigrid}
		S. Grepstad, G. Larcher, {\it On Pair Correlation and Discrepancy,} Arch. Math. (Basel) 109, no. 2, 143\,--\,149 (2017).
		
		\bibitem{our_other_paper}
		M. Hauke, A. Zafeiropoulos, {\it Poissonian Correlations of higher orders,}  Journal of Number Theory, Vol. 243, 202\,--\,240 (2023). 
		
		\bibitem{heath}
		D.R. Heath-Brown, {\it Pair correlation for fractional parts of $\alpha n^2$},
		Math. Proc. Cambridge Philos. Soc. 148, no.3, 385\,--\,407, (2010).
		
		\bibitem{kuipers}
		L. Kuipers, H. Niederreiter, {\it Uniform Distribution of Sequences,}
		Wiley (1974).
		
		\bibitem{kr}
		P. Kurlberg, Z. Rudnick, {\it The distribution of spacings between quadratic residues, } Duke J. Math., 100: 211\,--\,242 (1999).
		
		\bibitem{larch} 
		G. Larcher, W. Stockinger, {\it Some negative results related to Poissonian pair correlation problems,} Discrete Math. 343, no. 2, Article ID 111656, 11 p. (2020).
		
		\bibitem{lutsko}
		C. Lutsko,  N. Technau, {\it Full Poissonian Local Statistics of Slowly Growing Sequences.} Pre-print: arXiv 2206.07809.
		
		\bibitem{marklof}
		J. Marklof, \emph{ Pair correlation and equidistribution on manifolds,} Monatsh. Math. 191, 279\,--\,294 (2020).
		
		\bibitem{mark_stro} 
		J. Marklof, A. Str\"ombergsson, {\it Equidistribution of Kronecker sequences along closed horocycles.} Geom. Funct. Anal. 13, no.6, 1239\,--\,1280 (2003).
		
		\bibitem{marklof2}
		J. Marklof, N. Yesha, {\it Pair correlation for quadratic polynomials mod 1}, Compositio Mathematica 154, 960\,--\,983 (2018).
		
		\bibitem{tech_walk} 
		N. Technau, A. Walker, {\it On the triple correlations of fractional parts of $n^2 \alpha$}, Canad. J. Math. 74, no.5, 1347\,--\,1384, (2022).
		
		\bibitem{rud_sar} 
		Z. Rudnick, P. Sarnak,{\it The pair correlation function of fractional parts of polynomials}, Comm. Math. Phys. 194 (1998), no. 1, 61\,--\,70.  
		
		\bibitem{rudnick2}
		Z. Rudnick, P. Sarnak, A. Zaharescu, {\it The distribution of spacings between the fractional parts of $n^2\alpha$}, Invent. Math. 145, no. 1, 37\,--\,57 (2001).
		
		\bibitem{rudnick3}
		Z. Rudnick, A. Zaharescu, {\it  The distribution of spacings between fractional parts of lacunary sequences}, Forum Math. 14, 691\,--\,712, (2002).
		
		\bibitem{steinerberger}
		S. Steinerberger, {\it Poissonian Pair Correlation in Higher Dimensions,} J. Number Theory, 208: 47\,--\,58, (2020).  
		
		\bibitem{weiss}
		C. Wei{\ss}, {\it Some connections between discrepancy, finite gap properties, and pair correlations,} Monatsh. Math. 199, no. 4, 909\,--\,927 (2022). 
	\end{thebibliography}
\end{document}